\documentclass[journal,transmag,onecolumn]{IEEEtran}
\ifCLASSINFOpdf
  \usepackage[pdftex]{graphicx}
  \DeclareGraphicsExtensions{.pdf,.jpeg,.png}
\else
  \usepackage[dvips]{graphicx}
  \DeclareGraphicsExtensions{.eps}
\fi
\usepackage[cmex10]{amsmath}
\usepackage{amssymb} 
\usepackage{color}
\usepackage{tikz}	
\usepackage[european]{circuitikz} 
\allowdisplaybreaks

\newcommand{\eq}[1]{\eqref{eq:#1}}
\newcommand{\eqlabel}[1]{\label{eq:#1}}
\newcommand{\Fig}[1]{Figure~\ref{fig:#1}}
\newcommand{\fig}[1]{Fig.~\ref{fig:#1}}
\newcommand{\figlabel}[1]{\label{fig:#1}}
\newcommand{\Tab}[1]{Table~\ref{tab:#1}}
\newcommand{\tab}[1]{Tab.~\ref{tab:#1}}

\newcommand{\+}[2]{\def#1{#2}}
\newcommand{\1}[2]{\def#1##1{#2}}

\+{\ms}{\textrm{ms}}        
\+{\mV}{\textrm{mV}}        
\+{\us}{\mu\textrm{s}}      

\+{\uA}{\mu\textrm{A}}      
\+{\uF}{\mu\textrm{F}}      
\+{\cm}{\mathrm{cm}}        

\1{\card}{\#\,#1}           
\newcommand{\comm}[2]{{\left[{#1},{#2}\right]}}
\1{\D}{\times 10^{#1}}      
\+{\d}{\mathrm{d}}          
\1{\ddt}{\frac{\d{#1}}{\d{\t}}}             
\+{\e}{\mathrm{e}}          
\newcommand{\mat}[1]{\boldsymbol{#1}}       
\def\O(#1){\mathcal{O}\left(#1\right)}      
\newcommand{\T}{^{\!\top}} 

\providecommand{\abs}[1]{{\lvert}#1{\rvert}}
\providecommand{\norm}[1]{{\lVert}#1{\rVert}}

\1{\alp}{\alpha_{#1}}       
\1{\alY}{\alpha_{#1}(\Vm)}  
\1{\btY}{\beta_{#1}(\Vm)}   
\+{\alm}{\alpha_{\mgate}}   
\+{\btm}{\beta_{\mgate}}    
\+{\alh}{\alpha_{\hgate}}   
\+{\bth}{\beta_{\hgate}}    
\+{\Cm}{C}                  
\+{\dt}{\Delta t}           
\+{\dV}{\Delta V}           
\+{\E}{E}                   
\+{\Err}{\mathcal{E}}       
\+{\ErrFE}{\Err_{\mathrm{FE}}}   
\+{\ErrMRL}{\Err_{\mathrm{MRL}}} 
\+{\ErrHOS}{\Err_{\mathrm{HOS}}} 
\+{\ErrOS}{\Err_{\mathrm{OS}}}   
\+{\G}{G}                   
\+{\g}{g}                   
\+{\hgate}{h}               
\+{\I}{I}                   
\+{\Istim}{I_\mathrm{stim}} 
\+{\INa}{I_{\mathrm{Na}}}   
\+{\IKs}{I_{\mathrm{Ks}}}   
\+{\i}{i}                   
\+{\iset}{\mathcal{I}}      
\+{\j}{j}                   
\+{\jset}{\mathcal{J}}      
\+{\k}{k}                   
\+{\knum}{K}                
\+{\kset}{\mathcal{K}}      
\+{\knb}{\mathcal{K'}}      
\+{\krD}{\delta}            
\+{\l}{\ell}                
\+{\lnum}{L}                
\1{\M}{\mT{#1}}             
\+{\MEF}{\M{\mathrm{EF}}}   
\+{\MHOS}{\M{\mathrm{HOS}}} 
\+{\MMRL}{\M{\mathrm{MRL}}} 
\+{\MOS}{\M{\mathrm{OS}}}   
\+{\m}{m}                   
\+{\mA}{\mat{A}}            
\+{\matA}{\mat{X}}          
\+{\matB}{\mat{Y}}          
\+{\mB}{\mat{B}}            
\+{\mE}{\mat E}             
\+{\mD}{\mat{\Lambda}}      
\+{\mS}{\mat{S}}            
\1{\mT}{\mat{T}_{#1}}       
\+{\mgate}{m}               
\+{\n}{n}                   
\+{\Oref}{\sO^*}            
\1{\Pop}{P_{#1}}            
\+{\t}{t}                   
\+{\tf}{\t'}                
\+{\tmax}{\t_{\max}}        
\+{\tmin}{\t_{\min}}        
\1{\tC}{\tau_{#1}}          
\1{\tVm}{\tilde{V}_{#1}}    
\+{\tVmin}{V_{\min}}        
\+{\tVmax}{V_{\max}}        
\+{\u}{u}                   
\+{\Vm}{V_m}                
\1{\vu}{\vec{u}_{#1}}       
\+{\vX}{\vec{X}}            
\1{\y}{y_{#1}}              
\1{\yss}{\bar{y}_{#1}}      

\+{\C}{\mathrm{C}} 
\+{\IF}{\mathrm{IF}} 
\+{\IC}{\mathrm{IC}} 
\+{\IM}{\mathrm{IM}} 
\+{\CI}{\mathrm{CI}} 
\+{\OI}{\mathrm{OI}} 
\+{\RI}{\mathrm{RI}} 

\+{\sO}{O}
\+{\sP}{P}
\+{\sQ}{Q}
\+{\sR}{R}
\+{\sS}{S}
\+{\sT}{T}
\+{\sU}{U}
\+{\sV}{V}
\+{\sW}{W}

\+{\alRQ}{\alp{RQ}}
\+{\alQP}{\alp{QP}}
\+{\alPO}{\alp{PO}}
\+{\alST}{\alp{ST}}
\+{\alTU}{\alp{TU}}
\+{\alQR}{\alp{QR}}
\+{\alPQ}{\alp{PQ}}
\+{\alOP}{\alp{OP}}
\+{\alTS}{\alp{TS}}
\+{\alUT}{\alp{UT}}
\+{\alUP}{\alp{UP}}
\+{\alTQ}{\alp{TQ}}
\+{\alSR}{\alp{SR}}
\+{\alPU}{\alp{PU}}
\+{\alQT}{\alp{QT}}
\+{\alRS}{\alp{RS}}
\+{\alOU}{\alp{OU}}
\+{\alUO}{\alp{UO}}
\+{\alUV}{\alp{UV}}
\+{\alVU}{\alp{VU}}
\+{\alVW}{\alp{VW}}
\+{\alWV}{\alp{WV}}
\+{\alUN}{\alp{UN}}
\+{\alQN}{\alp{QN}}
\+{\alNU}{\alp{NU}}
\+{\alNQ}{\alp{NQ}}
\+{\alNO}{\alp{NO}}
\+{\alNP}{\alp{NP}}
\+{\alMS}{\alp{MS}}
\+{\alMT}{\alp{MT}}
\+{\alMU}{\alp{MU}}
\+{\alRM}{\alp{RM}}
\+{\alQM}{\alp{QM}}
\+{\alPM}{\alp{PM}}
\+{\alVM}{\alp{VM}}
\+{\alOM}{\alp{OM}}
\+{\alMO}{\alp{MO}}
\+{\alMP}{\alp{MP}}
\+{\alMQ}{\alp{MQ}}
\+{\alMR}{\alp{MR}}
\+{\alMV}{\alp{MV}}
\+{\alLR}{\alp{LR}}
\+{\alLQ}{\alp{LQ}}
\+{\alPL}{\alp{PL}}
\+{\alLP}{\alp{LP}}
\+{\alML}{\alp{ML}}
\+{\alLM}{\alp{LM}}

\+{\kPO}{K_{PO}}
\+{\kQO}{K_{QO}}
\+{\kRO}{K_{RO}}
\+{\kQP}{K_{QP}}
\+{\kRP}{K_{RP}}
\+{\kRQ}{K_{RQ}}
\+{\kST}{K_{ST}}
\+{\kSU}{K_{SU}}
\+{\kPU}{K_{PU}}
\+{\kQU}{K_{QU}}
\+{\kRU}{K_{RU}}

\+{\lOP}{L_{OP}}
\+{\lOQ}{L_{OQ}}
\+{\lOR}{L_{OR}}
\+{\lPQ}{L_{PQ}}
\+{\lPR}{L_{PR}}
\+{\lUS}{L_{US}}
\+{\lUT}{L_{UT}}

\+{\mult}{\mu}
\+{\mOP}{\mult_{OP}}
\+{\mOU}{\mult_{OU}}
\+{\mPO}{\mult_{PO}}
\+{\mPQ}{\mult_{PQ}}
\+{\mQP}{\mult_{QP}}
\+{\mQR}{\mult_{QR}}
\+{\mRQ}{\mult_{RQ}}
\+{\mST}{\mult_{ST}}
\+{\mTS}{\mult_{TS}}
\+{\mTU}{\mult_{TU}}
\+{\mUT}{\mult_{UT}}

\+{\orig}{\color{red}{-\#}}
\+{\update}{\color{cyan}{+\#}}
\+{\diff}{\hspace{1em}^\#}
\def\eqtwo(#1,#2){(\ref{eq:#1},\ref{eq:#2})}

\+{\Naion}            {\mathrm  {Na}}
\+{\Kion}             {\mathrm  {K}}
\+{\Caion}            {\mathrm  {Ca}}
\+{\NSR}            {\mathrm  {NSR}}
\+{\JSR}            {\mathrm  {JSR}}

\+{\MACacap}          {A_{Cap}}
\+{\MACacos}          {\arccos}
\+{\MACageo}          {A_{Geo}}
\+{\MACakOne}         {\alpha_{\Kion 1}}
\+{\MACb}             {b}
\+{\MACB}             {B}
\+{\MACbjsr}          {b_       {\JSR}}
\+{\MACbkOne}         {\beta_       {\Kion 1}}
\+{\MACbss}           {b_       {\infty}}
\+{\MACcai}           {\MACcamyo}
\+{\MACcaion}         {[\Caion^{2+}] _ {ion}}
\+{\MACcajsr}         {[\Caion^{2+}] _ {\JSR}}
\+{\MACcamyo}         {[\Caion^{2+}] _ {i}}
\+{\MACcamyot}        {[\Caion^{2+}] _ {myot}}
\+{\MACcansr}         {[\Caion^{2+}] _ {\NSR}}
\+{\MACcao}           {\MACcaout}
\+{\MACcaout}         {[\Caion^{2+}] _ {o}}
\+{\MACC}             {C}
\+{\MACcjsr}          {c_       {\JSR}}
\+{\MACCMDN}          {\mathrm{CMDN}}
\+{\MACcos}           {\cos}
\+{\MACCSQN}          {\mathrm{CSQN}}
\+{\MACCyclelength}   {\mathrm{CL}}
\+{\MACdcai}          {\Delta \MACcai }
\+{\MACdcajsr}        {\Delta \MACcajsr}
\+{\MACdcansr}        {\Delta \MACcansr}
\+{\MACd}             {d}
\+{\MACD}             {D}
\+{\MACdKin}          {\Delta[\Kion^+]_{i}}
\+{\MACdnain}         {\Delta[\Naion^+]_{i}}
\+{\MACdss}           {d_{\infty}}
\+{\MACdt}            {\Delta t}
\+{\MACdvdt}          {\ddt{\MACv}}
\+{\MACdvdtmax}       {\ddt{\MACv}_{max}}
\+{\MACdvolderdt}     {\d \MACv_{olderdt}}
\+{\MACdvoldestdt}    {\d \MACv_{oldestdt}}
\+{\MACeca}           {E_{\Caion }}
\+{\MACEca}           {E_{\Caion }}
\+{\MACEkOne}         {E_{\Kion 1}}
\+{\MACekr}           {E_{\Kion r}}
\+{\MACEks}           {E_{\Kion s}}
\+{\MACEna}           {E_{\Naion }}
\+{\MACexp}           {\exp}
\+{\MACFab}           {F_{ab}}
\+{\MACfca}           {f_{\Caion }}
\+{\MACf}             {f}
\+{\MACF}             {F}
\+{\MACfNaK}          {f_{\Naion \Kion}}
\+{\MACfrdy}          {F}
\+{\MACfss}           {f_{\infty}}
\+{\MACgacai}         {\gamma_{\Caion i}}
\+{\MACgacao}         {\gamma_{\Caion o}}
\+{\MACgaki}          {\gamma_{\Kion i}}
\+{\MACgako}          {\gamma_{\Kion o}}
\+{\MACgammas}        {\eta}                            
\+{\MACganai}         {\gamma_{\Naion i}}
\+{\MACganao}         {\gamma_{\Naion o}}
\+{\MACgcat}          {\bar{G}_{\Caion (T)}}
\+{\MACg}             {g}
\+{\MACggs}           {g_{\infty}}
\+{\MACGkOne}         {\bar{G}_{\Kion 1}}
\+{\MACgkr}           {\bar{G}_{\Kion r}}
\+{\MACGks}           {\bar{G}_{\Kion s}}
\+{\MACGNa}           {G_{\Naion }}
\+{\MACGrel}          {G_{rel}}
\+{\MACibarca}        {\bar{I}_{\Caion }}
\+{\MACibark}         {\bar{I}_{\Caion \Kion }}
\+{\MACibarna}        {\bar{I}_{\Caion \Naion }}
\+{\MACicab}          {I_{\Caion b}}
\+{\MACicaL}          {I_{\Caion L}}
\+{\MACicaT}          {I_{\Caion (T)}}
\+{\MACIkOne}         {I_{\Kion 1}}
\+{\MACik}            {I_{\Kion }}
\+{\MACikns}          {I_{ns\Kion}}
\+{\MACIkp}           {I_{\Kion p}}
\+{\MACikr}           {I_{\Kion r}}
\+{\MACiks}           {I_{\Kion s}}
\+{\MACilca}          {I_{\Caion }}
\+{\MACilcak}         {I_{\Caion \Kion}}
\+{\MACilcana}        {I_{\Caion \Naion}}
\+{\MACIleak}         {I_{leak}}
\+{\MACinab}          {I_{\Naion b}}
\+{\MACInaca}         {I_{\Naion \Caion }}
\+{\MACina}           {I_{\Naion }}
\+{\MACInak}          {I_{\Naion \Kion}}
\+{\MACinsca}         {I_{ns(\Caion) }}
\+{\MACInsk}          {\bar{I}_{ns\Kion }}
\+{\MACinsna}         {I_{ns\Naion }}
\+{\MACInsna}         {\bar{I}_{ns\Naion }}
\+{\MACipca}          {I_{p(\Caion) }}
\+{\MACIrel}          {I_{rel}}
\+{\MACitca}          {I_{t\Caion }}
\+{\MACit}            {I_{t}}
\+{\MACitk}           {I_{t\Kion }}
\+{\MACitna}          {I_{t\Naion }}
\+{\MACito}           {I_{to}}
\+{\MACItr}           {I_{tr}}
\+{\MACIup}           {I_{up}}
\+{\MACiv}            {I_{v}}
\+{\MACkOne}          {\Kion1_{\infty}}
\+{\MACki}            {\MACKin}
\+{\MACKin}           {[\Kion^+]_{i}}
\+{\MACKleak}         {K_{ leak}}
\+{\MACkmca}          {K_{ m\Caion} }
\+{\MACko}            {\MACKout}
\+{\MACKout}          {[\Kion^+]_{o}}
\+{\MACkp}            {K_{ p}}
\+{\MACl}             {L}
\+{\MAClog}           {\log}
\+{\MACnai}           {\MACNain}
\+{\MACNain}          {[\Naion^+]_{i}}
\+{\MACnao}           {\MACNaout}
\+{\MACNaout}         {[\Naion^+]_{o}}
\+{\MACNumAPs}        {NumAPs}
\+{\MACpca}           {P_{\Caion }}
\+{\MACpi}            {\pi}
\+{\MACpk}            {P_{\Kion }}
\+{\MACpna}           {P_{\Naion }}
\+{\MACPnsca}         {P_{ns(\Caion) }}
\+{\MACprnak}         {P_{\Naion \Kion}}
\+{\MACrad}           {r}
\+{\MACRestinterval}  {\mathrm{Rest}_{interval}}
\+{\MACrot}           {R_{\Kion r}}
\+{\MACr}             {R}
\+{\MACR}             {\MACr}
\+{\MACryrclose}      {\mathrm{ryr}_{close}}
\+{\MACryropen}       {\mathrm{ryr}_{open}}
\+{\MACsig}           {\sigma}
\+{\MACsqrt}          {\mathrm{sqrt}}
\+{\MACstimulus}      {I_{st}}
\+{\MACtaub}          {\tau_{b}  }
\+{\MACtaud}          {\tau_{d}  }
\+{\MACtauf}          {\tau_{f}  }
\+{\MACtaug}          {\tau_{g}  }
\+{\MACtauxr}         {\tau_{xr} }
\+{\MACtauxsOne}      {\tau_{xs1}}
\+{\MACtauxsTwo}      {\tau_{xs2}}
\+{\MACtc}            {t_{c}}
\+{\MACtemp}          {\MACT}
\+{\MACTimeofdvdtmax} {\mathrm{Time_{ofdvdtmax}}}
\+{\MACTRPN}          {\mathrm{TRPN}}
\+{\MACt}             {t}
\+{\MACT}             {T}
\+{\MACvcell}         {V_{cell}}
\+{\MACvjsr}          {V_{\JSR}}
\+{\MACvmyo}          {V_{myo}} 
\+{\MACvnsr}          {V_{\NSR}}
\+{\MACvolder}        {\MACv_{older}}
\+{\MACvoldest}       {\MACv_{oldest}}
\+{\MACvold}          {\MACv_{old}}
\+{\MACv}             {\Vm}
\+{\MACxrss}          {X_{r\infty}}
\+{\MACxr}            {X_{r}}
\+{\MACxsOness}       {x_{s1\infty}}
\+{\MACxsOne}         {x_{s1}}
\+{\MACxsTwoss}       {x_{s2\infty}}
\+{\MACxsTwo}         {x_{s2}}
\+{\MACzca}           {z_{\Caion }}
\+{\MACzk}            {z_{\Kion }}
\+{\MACzna}           {z_{\Naion }}

\begin{document}
\title{Exponential integrators for a Markov chain model \\ of the fast sodium channel of cardiomyocytes}
\author{%
  \IEEEauthorblockN{
    Tom\' a\v s Star\' y\IEEEauthorrefmark{1}, %
    Vadim N. Biktashev\IEEEauthorrefmark{1} %
  }\IEEEauthorblockA{\IEEEauthorrefmark{1}%
    College of Engineering, Mathematics and Physical Sciences, %
    University of Exeter, %
    EX4 4QF, UK %
  }\thanks{
    Manuscript received ??; %
    revised ??. %
    Corresponding author: \sV. N. Biktashev
  }
}

\markboth{Paper submitted to IEEE Trans Biomed Eng}%
{Stary and Biktashev: Exponential integrators for a Markov chain model \ldots}

\IEEEtitleabstractindextext{%
\begin{abstract}
  The modern Markov chain models of ionic channels in excitable
  membranes are numerically stiff. The popular numerical methods for
  these models require very small time steps to ensure
  stability. Our objective is to formulate and test two methods
  addressing this issue, so that the timestep can be chosen
  based on accuracy rather than stability.
  Both proposed methods extend Rush-Larsen technique, which was
  originally developed to Hogdkin-Huxley type gate models.  One
  method, ``Matrix Rush-Larsen'' (MRL) uses a matrix reformulation
  of the Rush-Larsen scheme, where the matrix exponentials are
  calculated using precomputed tables of eigenvalues and
  eigenvectors. The other, ``hybrid operator splitting'' (HOS) method 
  exploits asymptotic properties of a particular Markov chain
  model, allowing explicit analytical expressions for the substeps.
  We test both methods on the Clancy and Rudy (2002) $\INa$ Markov chain
  model. With precomputed tables for functions of the transmembrane
  voltage, both methods are comparable to the forward Euler method in
  accuracy and computational cost, but allow longer time steps without
  numerical instability.
  We conclude that both methods are of practical interest. MRL
  requires more computations than HOS, but is formulated
  in general terms which can be readily extended to other Markov Chain
  channel models, whereas the utility of HOS depends on the asymptotic
  properties of a particular model.
  The significance of the methods is that they allow a considerable
  speed-up of large-scale computations of cardiac excitation models by
  increasing the time step, while maintaining acceptable accuracy and
  preserving numerical stability.
\end{abstract}

}

\maketitle

\IEEEdisplaynontitleabstractindextext
\IEEEpeerreviewmaketitle

\section{Introduction}

\IEEEPARstart{M}{athematical} models are an essential part of the
modern cardiac electrophysiology. They are used for hypothesis testing
in research and as a guide for clinical decision. A typical definition
of such a model is a high-dimensional (tens of equations) system of
ordinary differential equations per excitable unit. Detailed
simulations of the heart involve solving such systems for each of
millions of cells placed in a mesh representing the cardiac
tissue. Such large-scale models can be computationally extremely
expensive, hence significant efforts are directed to develop efficient
numerical methods for solving such systems.

A typical cardiac excitation model is centered around the 
Kirchhoff circuit law which gives
\begin{equation}
  \Istim(\t) = \Cm \ddt{\Vm} + \sum_\l \I_\l,  \eqlabel{current-law}
\end{equation}
where $\Cm$ is the cell membrane's capacitance, $\Vm=\Vm(\t)$ is the
transmembrane potential difference, and $\I_\l$, $\l=1,\ldots,\lnum$, are
currents through ion-specific channels. The currents, in turn, are
determined by the Ohm's law,
\begin{align}
  \I_\l = \G_\l \Pop{\l}(\t) \left[ \Vm(\t) - \E_\l(\vX(\t)) \right]
\end{align}
where $\E_\l(\vX)$ is the ion-specific electromotive force, depending
on the ionic concentrations $\vX$ via the Nerst equation, $\G_\l$ is
the total conductance of channels of type $\l$ when they
are all open, and $\Pop{\l}$ is the probability of those channels to
be open.

The components of the vector $\vX$ are intra- and extra-cellular ionic
concentrations, which change in time in the obvious way in accordance
with the ionic fluxes and the corresponding volumes; some
concentrations in some models are assumed constant. The dynamics of
the open probabilities is much more nontrivial, as it reflects the
conformation dynamics of the proteins constituting the ion channels .

The classical description of these dynamics, going back to Hodgkin and
Huxley (1952)~\cite{Hodgkin-Huxley-1952}, has the form
\begin{align}
  \Pop{\l}(\t) = \prod_{\i\in\iset(\l)} \y{\i}
\end{align}
with a popular, although different from the original Hodgkin and
Huxley's, interpretation that the set $\iset(\l)$ corresponds the
subunits of the channel, called ``gates''. These subunits are assumed
statistically independent, each of them can be either in an ``open''
or a ``closed'' state, and the channel is open if and only if each of
the subunits is open. Variables $\y{i}$ then are open probabilities of
the gates, and their dynamics are described by
\begin{align}
  \eqlabel{gate}
  \ddt{\y{\i}} = \alY{\i} (1-\y{\i}) - \btY{\i} \y{\i} ,
\end{align}
where $\alY{\i}$ are opening rates and $\btY{\i}$ are closing rates. 

For instance, the original Hodgkin-Huxley description of the fast
sodium current ($\INa$) channel uses $\card{\iset(\INa)}=4$ gates, three
of which, called m-gates, have identical opening
$\alY{1}=\alY{2}=\alY{3}=\alm(\Vm)$ and closing
$\btY{1}=\btY{2}=\btY{3}=\btm(\Vm)$ rates, and the fourth, called
h-gate, has rates $\alY{4}=\alh(\Vm)$ and $\btY{4}=\bth(\Vm)$,
hence for this case we have 
\begin{align}
  & \Pop{\INa}(\t) = \mgate^3 \hgate, \\
  & \ddt{\mgate} = \alm(\Vm) (1-\mgate) - \btm(\Vm) \mgate, \\
  & \ddt{\hgate} = \alh(\Vm) (1-\hgate) - \bth(\Vm) \hgate.
\end{align}

A more recent approach is modelling the probabilities of the channel
molecules, as a whole, to be in specific conformation states, without
the restricting assumptions of statistically independent subunits and
only two states for any subunit. This gives generic continuous time
Markov chain (MC) models. Let the probability of the $\l$'th channel to be in
the $\k$'th state at time $\t$ be $\u_\k(\t)$ (``state occupancy''),
$\k\in\kset(\l)$, $\knum(\l)=\card{\kset(\l)}$, and all such
probabilities be considered components of the state (column-) vector
$\vu{}=\left(\u_\k\right)\T=\vu{}(\t)$. Let $\g_\k$ be the relative permeability
of the state $\k$, then we have
\begin{align}
  \eqlabel{MC-Pop}
  \Pop{\l}(\t) = \sum\limits_{\k=\kset(\l)} \g_\k \u_\k(\t) .
\end{align}
Typically, $\g_\k=\krD_{\k,\k^*}$ where $\k^*$ is the ``open state''.
The time evolution of the state vector is described by the
system of linear ODEs,
known in particular as \emph{Kolmogorov (forward) equations}, or
\emph{master equation}, of the form
\begin{align}
  \eqlabel{markov-generic}
	\ddt{\vu{}} = \mA(\Vm) \vu{} ,
\end{align}
where the non-diagonal components of the matrix $\mA(\Vm)$ are the
transition rates (TR) between the states, and the diagonal components are
defined by the condition $\sum_{\k\in\kset}\u_\k=1$ and consequently
sum of any column of $\mA$ should be zero. 

The ODE system for cellular membrane can be solved on a computer using
standard numerical solvers. A typical solver
iteratively computes the states of the system using
time-stepping algorithms, that is computing the states 
at times $\t_\n = \t_0 + \n\dt$. The size of the time step
$\dt$ is inversely proportional to the computational cost, measured
as CPU time required for the computation. Increasing the time step is
a straightforward way of reducing the computational cost. 

The maximal acceptable time step is limited by considerations of
accuracy and stability (see e.g.~\cite[Sections
  5.10, 5.11]{Burden-Faires-textbook}).  Whereas the former is ``relative''
in that it depends on the aims of the research, the latter has a more
``absolute'' character in that if stability conditions are not
satisfied, the solution is unusable for any purpose. Typically, when
the time step exceeds the stability limit, the numerical solution is
characterized by wild oscillations around the exact solution, and
quite often will lead to numerical overflow.

\begin{figure}[b!]
  \centering
  \includegraphics[width=3.5in]{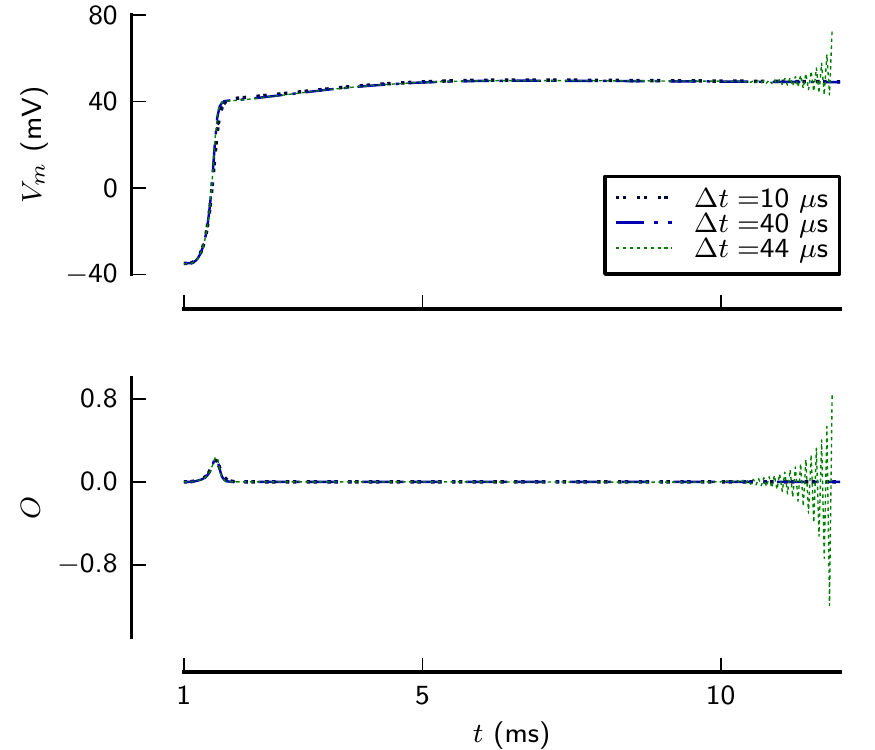}
  \caption[]{%
    Instability of $\INa$ Markov chain model at longer time steps. %
    The model was solved using forward Euler method using three
    different step sizes: $\dt=10\,\us$, $\dt=40\,\us$, and
    $\dt=44\,\us$. %
    The top panel shows the membrane potential ($\Vm$), %
    the bottom panel shows state occupancy of the open state $\sO$. %
  }
  \figlabel{stiffness}
\end{figure}

 Simple explicit solvers suffer from instabilities the most,
  and implicit, stable methods, applicable to generic systems of ODEs,
  are complicated and often costly. The motivation for our research
was that taking into account the specific properties of the problem
can offer some advantages. Specifically, we have in mind two distinct
considerations.

One consideration is that the TRs
can range through several orders of magnitude, and some of them can be
much faster than other processes described by the excitable
cell model. This split of the speeds of the variables suggest a
possibility to exploit asymptotic methods.

The other consideration is the linearity of the
system~\eq{markov-generic}. Here we are inspired by the example of the
exponential integrator algorithm developed by Rush and Larsen in
1978~\cite{Rush-Larsen-1978}. It is based on the assumption that the
transmembrane voltage, on which the TRs in the gate model \eq{gate}
depend, changes only slightly during one time step $\dt$. So during
one time step, the TRs can be approximated by constants, and
the equation~\eq{gate} is then solved analytically. The solution can
be conveniently defined in terms of the ``steady state'' 
$\yss{\i} = \alY{\i}/[\alY{\i}+\btY{\i}]$ 
and the ``time constant'' 
$\tC{\i} = 1/[\alY{\i}+\btY{\i}]$ at a
given potential $\Vm$ presumed constant:
\begin{align}
  \y{\i}^{\n+1} = \yss{\i}(\Vm) - [\yss{\i}(\Vm) - \y{\i}^\n]
  \exp\left(-\frac{\dt}{\tC{\i}(\Vm)}\right) .
\end{align}

The Rush-Larsen (RL) scheme is easy to implement, gives good results
and is very popular in computational cardiac
electrophysiology. Its stability and approximation properties
  have been extensively discussed in literature,
  including its relation to general exponential integrators family, its
  extension beyond gating variables by linearization, and improving
  its approximation properties, see e.g.~\cite{%
    Sundnes-etal-2009,%
    Perego-Veneziani-2009,%
    Marsh-etal-2012}.  However it is designed for a single ODE and is
not immediately applicable for  MC models which are systems of
coupled ODEs. And yet MC models are known to suffer from severe numerical
instability issues, just as, or even more than, the gate models
(\fig{stiffness}). %
  The classical techniques for numerical solution of continuous-time
  MC models involve finding the eigenvalues and associated
  eigenvectors of the transition matrix. Direct implementation of this
  approach to very large MCs is problematic, see
  e.g.~\cite{Reidman-Trivedi-1988}. However the MCs describing ionic
  channels are relatively small so the direct approach is feasible.

In this paper, we discuss two methods for numerical solution of MC
models based on these two considerations.

\section{Methods}

\subsection{Models}

To test the suggested numerical methods we have chosen the MC model of
the $\INa$ channel  by Clancy and Rudy~\cite{Clancy-Rudy-2002}
(\fig{diag_MC}), which is one of the most popular MC models.
We used the formulation of the MC model and the whole cell model into
which it was incorporated, as implemented in the authors' code kindly
provided by C.E.~Clancy. It most
closely corresponds to the Luo-Rudy model~\cite{Luo-Rudy-1994a} with
modifications described in~\cite{Zeng-etal-1995,Viswanathan-etal-1999a}, and some
further minor differences. For the sake of reproducibility of our
results, we describe the whole model  in the supplementary
material, highlighting all the differences from the published models
that we have detected. For the same purpose, we put a simplified
version of the C code we used in the simulations described below in
the supplementary materials.

\begin{figure}
  \centering
    \includegraphics{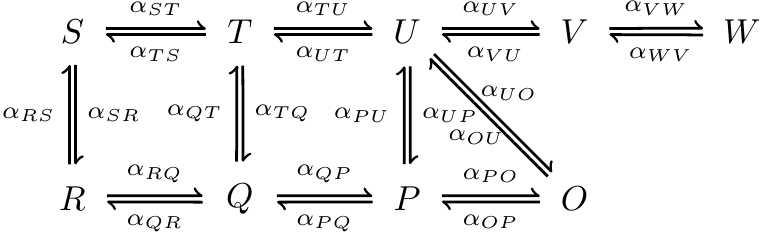}
  \caption[]{Markov chain model of $\INa$ channel.}
  \figlabel{diag_MC}
\end{figure}

\begin{table}
  \caption[]{Change of states variables terminology.}
  \centering
    \begin{tabular}{lll}
      \hline 
      Standard  & Our & initial value \\      
      \hline
      $\sO$     &    $\sO$ &$4.386\D{-8}$	 \\ 
      $\C_1$  &    $\sP$ &$5.329\D{-5}$	 \\
      $\C_2$  &    $\sQ$ &$1.064\D{-2}$	 \\
      $\C_3$  &    $\sR$ &$8.018\D{-1}$	 \\
      $\IC_3$ &    $\sS$ &$1.436\D{-1}$	 \\
      $\IC_2$ &    $\sT$ &$1.907\D{-3}$ 	 \\
      $\IF$   &    $\sU$ &$1.111\D{-5}$	 \\
      $\IM_1$ &    $\sV$ &$8.417\D{-4}$	 \\
      $\IM_2$ &    $\sW$ &$4.118\D{-2}$     \\
      \hline
    \end{tabular} \label{tab:states-term}
\end{table}

For convenience, we changed the notation for the MC states and TRs. The
states were named in alphabetical order, starting with $\sO$ for the
open state, in a clockwise direction as in
\fig{diag_MC}. See~\tab{states-term} and \tab{rates-term} for the correspondence with the
original notation.
The model contains 9 interconnected states. The state $\sO$ represents the
conformation of the ion channel that allows the flow of ions between
the intracellular and extracellular environment. The remaining states
($\sP$, $\sQ$, $\sR$, $\sS$, $\sT$, $\sU$, $\sV$ and $\sW$) represent non-conductive
conformations of the channel, so we can say that for this model
$\g_\k=\krD_{k,1}$, where $\u_1=\sO$. 
There are 11
possible bidirectional transitions between states, but some of the
corresponding 22 TRs are described by identical functions, so there
are only 14 distinct TR definitions.
We denote the TRs by $\alp{}$ with a subscript showing the direction of the
  transition, e.g. $\alPO$ is the transition rate from state $\sP$
  into state $\sO$. See~\tab{rates-term} for the link with the
  original notations.

\begin{table}
  \caption[]{Change of transition rates (TR) terminology.}
  \centering
    \begin{tabular}{ll}
      \hline 
      Standard  & Our  \\      
      \hline
      $\alpha_{11}$&      $\alRQ$, $\alST$                \\
      $\alpha_{12}$&      $\alQP$, $\alTU$                \\
      $\alpha_{13}$&      $\alPO$                              \\
      $\beta_{11}$ &      $\alQR$, $\alTS$                \\
      $\beta_{12}$ &      $\alPQ$, $\alUT$                \\
      $\beta_{13}$ &      $\alOP$                              \\
      $\alpha_2$   &      $\alOU$                              \\
      $\beta_2$    &      $\alUO$                              \\
      $\alpha_3$   &      $\alUP$, $\alTQ$, $\alSR$ \\
      $\beta_3$    &      $\alPU$, $\alQT$, $\alRS$ \\
      $\alpha_4$   &      $\alUV$                              \\
      $\beta_4$    &      $\alVU$                              \\
      $\alpha_5$   &      $\alVW$                              \\
      $\beta_5$    &      $\alWV$                              \\
      \hline
    \end{tabular} \label{tab:rates-term}
\end{table}

The TRs are shown on \fig{TR} as functions of the transmembrane potential
$\Vm$ in a physiologically relevant range. The values of TRs vary
across several orders of magnitude, from $10^{-11}\,\ms^{-1}$
to $10^2\,\ms^{-1}$. Some of the TRs are high at the lower
potentials, some are fast at higher potentials, and some are uniformly
low.

\begin{figure*}
  \centering
  \includegraphics[width=7.16in]{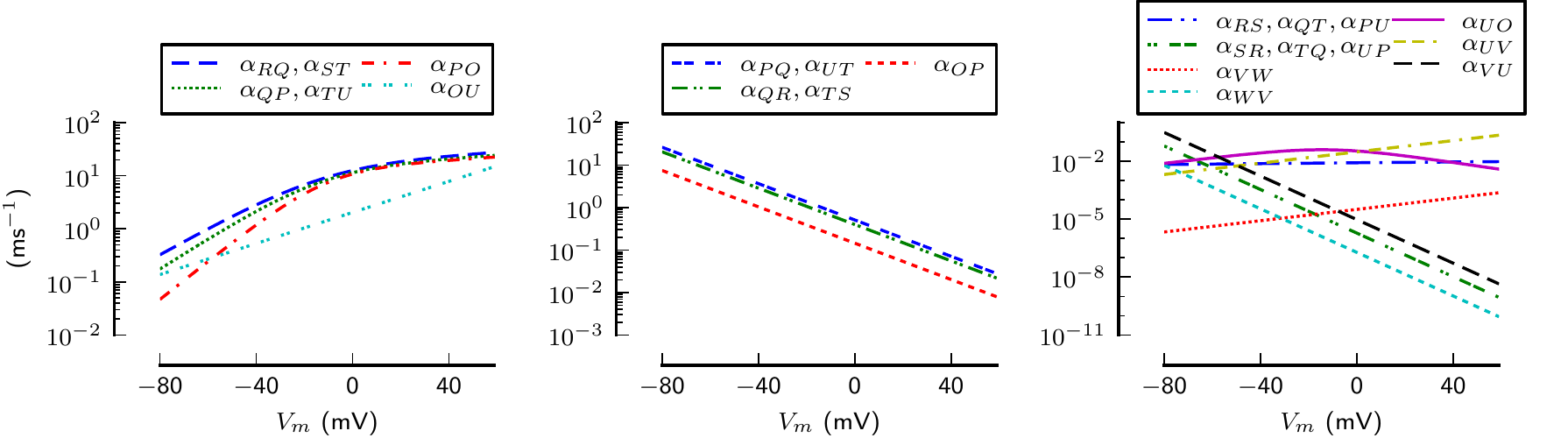}
  \caption[]{Transition rates (TR) of $\INa$ Markov chain model. %
    The left panel shows fast TR at high potentials: subsystem $\mA_0$;  %
    the middle panel shows fast TR at low potentials: subsystem $\mA_1$; %
    and the right panel shows slow TR in the entire range of
    potentials: subsystem $\mA_2$. %
  }
  \figlabel{TR}
\end{figure*}

The conductive (open) state $\sO$ is the only state that has immediate
effect on the $\INa$ current. 
The remaining 8 states of the model can
affect the current only indirectly by transitions to the open state
$\sO$. The time evolution of a generic state occupancy state $\u_\k$ is
described by a differential equation of the form
\begin{align}
  \ddt{\u_\k} = \sum\limits_{\k'\in\knb(\k)} 
  \left( \alp{\k',\k} \u_{\k'} - \alp{\k,\k'} \u_\k \right) ,
                                                  \eqlabel{genX}
\end{align}
where $\knb(\k)$ is the set of all the states interconnected with state
$\k$, which can be readily found from the diagram.  For example,
the occupancy of the open state, $\sO$, is described by the following
equation:
\begin{align*}
  \ddt{\sO} = \alPO\,\sP + \alUO\,\sU - (\alOP + \alOU)\,\sO . 
\end{align*}

By taking the sum of the equations \eq{genX} for all $\k$, one can see
that the sum of the right-hand sides equals to zero, and therefore the
system observes a states conservation law, which is consistent with
the definition of $\u_\k$ as probabilities, implying
$\sum_{\k\in\kset}\u_\k=1$. This is of course a generic property of a
continuous Markov chain. So the differential equations in the model
are not independent, which creates a possibility of reducing the
number of equations from 9 to 8, by computing one of the occupancies
through the conservation law rather than from its differential
equation. However the computational gain from this is insignificant,
and instead we used any deviations from the conservation law as an
indicator of the accuracy of the computations.

\subsection{Numerical Methods}

\subsubsection{Forward Euler}

The standard forward Euler (FE) method is the simplest timestepper for
differential equations. It defines the solution at the next time step,
$\vu{\n+1}=\vu{}(\t_{\n+1})$, in terms of the same at
the previous time step, $\vu{\n}=\vu{}(\t_{\n})$, using one-step
forward-time finite different approximation of the time derivative,
which for the system~\eq{markov-generic} gives
\begin{align}
  \vu{\n+1} = 
    \vu{\n} + \dt \mA(\Vm(\t_\n)) \vu{\n}.
\end{align}
The time discretization step $\dt = \t_{\n+1} - \t_\n$ is
presumed here the same for all steps of a simulation.

\subsubsection{Matrix Rush-Larsen}

The proposed Matrix Rush-Larsen method (MRL) assumes that the matrix
$\mA(\Vm)$ changes only slightly during one time step and
therefore can be approximated by a constant. The solution of \eq{markov-generic} can then be
written in terms of the matrix exponential,
\begin{align}
  \eqlabel{MRL-1}
  \vu{\n+1} = \exp\left[\mA(\Vm(\t_\n))\,\dt \right]\,\vu{\n} .
\end{align}
We assume that the matrix $\mA(\Vm)$ is diagonalizable, i.e. can be
represented in the form $\mA(\Vm) = \mS(\Vm) \mD(\Vm) \mS(\Vm)^{-1}$,
where matrix $\mS(\Vm)$ is composed of the eigenvectors concatenated
as column vectors, and matrix $\mD(\Vm)$ contains eigenvalues placed
on the corresponding places on the diagonal. %
  A sufficient condition of diagonalizability of a matrix is that all
  its eigenvalues are distinct, and this is the generic situation; %
but we of course check that it actually takes place in every
case. Then the matrix exponential is calculated as
\begin{align}
  \vu{\n+1} = \mS(\Vm)\,\exp\left( \mD(\Vm) \dt
  \right)\,\mS(\Vm)^{-1} \vu{\n}, 
\end{align}
where the exponential $\exp\left[\mD(\Vm) \dt\right]$ of the diagonal
matrix $\mD(\Vm)\dt$ is obtained by exponentiation of its diagonal
elements.

As the numerical solution of the eigenvalue problem is computationally
expensive, we precompute the matrices $\mS(\tVm{\j})$,
$\mS(\tVm{\j})^{-1}$ and $\mD(\tVm{j})$ for a fine grid of
physiological potentials, $\tVm{\j} = \tVmin + \j\dV$, $\j \in \jset =
\{0,1,\hdots \j_\mathrm{max} \}$, $\tVm{\j} \leq \tVmax$,
$\tVmin=-100$, $\tVmax=70$, $\dV=0.01$ (all in \mV) before compile
time and save them in a file. 

At start time, the eigenvalue and eigenvector matrices are loaded from
the file and we precompute, for $\dt$ used in the particular
simulation, the transition matrices
\begin{align}
  \mT{\j} = \mT{}\left(\tVm{\j}\right) = 
  \mS\left(\tVm{\j}\right)
  \exp\bigg[ \mD\left(\tVm{\j}\right) \dt \bigg]
  \mS\left(\tVm{\j}\right)^{-1}
\end{align}
for all $\j\in\jset$.
At the run time, the solver simply refers to the tabulated transition matrix $\mT{\j}$,
\begin{align}
  \vu{\n+1} = \mT{\j(n)} \vu{\n}
\end{align}
where $\tVm{\j(\n)}$ is the tabulated transmembrane potential that is the
nearest to $\Vm(\t_\n)$.

Along with the code, we provide precomputed files for a
  voltage step size of $\dV = 0.1$ (size of 4.85 MB), that are
  sufficient to obtain accurate results. The tables with $\dV =
  0.01\,\mV$ of 48.5 MB size, used for the simulations presented, are
  available from the authors upon request.

The method of tabulation (tab.) can be applied to all the presented
numerical methods.  However, its benefit is most essential in the MRL
method, because matrix exponentiation is computationally expensive. The
accuracy of the tabulation is dependent on the voltage step
(here 0.01\,\mV) which is a matter of choice depending on
memory availability and allowable pre-compile and start-time
computation time.

\subsubsection{Hybrid Operator Splitting}

The MRL achieves the purpose in principle but is relatively costly as
multiplication by a dense $\knum\times \knum$ matrix $\mT{\j(n)}$ is
required at each time step.  On the other hand, it did not at all
exploit the specific structure of the TR, illustrated by \fig{TR},
that is, that the matrix $\mA(\Vm)$ is sparse and some TRs are much
faster than others for some voltage ranges.  Hence we propose a 
hybrid operator splitting method (HOS), which combines FE and MRL,
and exploits the asymptotic structure of the TRs.  In this
  method, we set 
\begin{align*}
    \mA(\Vm) = \mA_0(\Vm)+\mA_1(\Vm)+\mA_2(\Vm)
\end{align*}
as described in \fig{TR}: $\mA_0$ contains
    only TRs that are fast at high values of $\Vm$ ($\alRQ$, $\alST$,
    $\alQP$, $\alTU$, $\alPO$ and $\alOU$); $\mA_1$ contains only TRs
    that are fast at low values of $\Vm$ ($\alPQ$, $\alUT$, $\alQR$,
    $\alTS$ and $\alOP$); and $\mA_2$ contains only uniformly slow TRs
    ($\alRS$, $\alQT$, $\alPU$, $\alSR$, $\alTQ$, $\alUP$, $\alVW$,
    $\alWV$, $\alUO$, $\alUV$ and $\alVU$). Explicit expressions for
    $\mA_j$ are given in the Supplement.

Every timestep is then done in three substeps,
\begin{align}
  \vu{\n+1/3} =&\, \exp(\dt \mA_0(\Vm(\t_\n))) \, \vu{\n} , \eqlabel{sol-subs-1} \\
  \vu{\n+2/3} =&\, \exp(\dt \mA_1(\Vm(\t_\n))) \, \vu{\n+1/3} , \eqlabel{sol-subs-2} \\
  \vu{\n+1} =&\, \vu{\n+2/3} +\dt \mA_2(\Vm(\t_\n)) \, \vu{\n+2/3} .\eqlabel{sol-subs-3}
\end{align}
In our case the matrix exponentials in the two fast subsystems
\eq{sol-subs-1} and \eq{sol-subs-2} are found analytically, through
solving the corresponding ODE systems.  This is possible because some
of the equations corresponding to the matrices
$\mA_0(\Vm)$ and $\mA_1(\Vm)$ are coupled in a specific manner and can
be solved one by one where solution of one equation is substituted in
the next etc. The full expressions and the method of
  derivation are given in
the Supplement; here we present the solution for state $\sO$ in the
equation \eq{sol-subs-1} as an example:
\begin{align*}
      \sO_{\n+1/3}= \mOU \sO_{\n} + \kPO \sP_{\n} + \kQO  \sQ_{\n} + \kRO  \sR_{\n},
\end{align*}
where 
\begin{align*}
  \kPO =&\,  \frac{\alPO(\mPO - \mOU)  }{\alOU - \alPO}, \\
  \kQO =&\, \frac{\alPO \alQP(\mQP - \mOU) }{(\alPO - \alQP)(\alOU - \alQP)} \nonumber \\
  &- \frac{\alPO \alQP (\mPO - \mOU)}{(\alPO - \alQP)(\alOU - \alPO)} , \\
  \kRO =&\,  - \frac{\alPO \alQP \alRQ(\mQP - \mOU) }{(\alQP - \alRQ)(\alPO - \alQP)(\alOU -\alQP)} 
  \nonumber \\ 
  &+ \frac{\alPO \alQP \alRQ(\mPO - \mOU) }{(\alQP - \alRQ)(\alPO -
    \alQP)(\alOU -\alPO)}      \nonumber \\
  & +\frac{\alPO \alQP \alRQ (\mRQ - \mOU)}{(\alQP - \alRQ)(\alPO -
    \alRQ)(\alOU - \alRQ)}  \nonumber \\
  &-\frac{\alPO \alQP \alRQ(\mPO - \mOU) }{(\alQP - \alRQ)(\alPO -
    \alRQ)(\alOU - \alPO)} , \nonumber \\
  \mult_{\j\k} &=\,\e^{ -\alp{\j\k} \dt }.
\end{align*}
%
The slow subsystem \eq{sol-subs-3} uses FE, and since it contains only
uniformly slow TRs, it can tolerate large time-steps, allowed by other
components of the cell model, without loss of stability.

\subsection{A priori error estimates}

Estimates by standard methods
(e.g.~\cite[Chapter 5]{Burden-Faires-textbook}, see details in the
Supplement), show that all three numerical schemes have local
truncation errors of the second order, i.e. $\Err\dt^2+\O(\dt^3)$,
although the coefficients $\Err$ vary: for FE we have %
$\ErrFE\leq\frac12\left(\norm{\mA}^2+\norm{\d\mA/\d\Vm}\abs{\d\Vm/\d\t}\right)$, %
for MRL we have %
$\ErrMRL=\frac12\norm{\d\mA/\d\Vm}\abs{\d\Vm/\d\t}$, %
and for HOS it is composed of contributions of the three substeps plus
the error due to operator splitting,
$\ErrHOS\leq
\frac12\abs{\d\Vm/\d\t}\left(
  \norm{\d\mA_0/\d\Vm}
 +\norm{\d\mA_1/\d\Vm}
 +\norm{\d\mA_2/\d\Vm}
\right)+\frac12||\mA_2||^2+\ErrOS$, %
where
$\ErrOS=\frac12\norm{\comm{\mA_1}{\mA_0}+\comm{\mA_2}{\mA_0}+\comm{\mA_2}{\mA_1}}$,
and $\comm{\matA}{\matB}\equiv\matA\matB-\matB \matA$.  %
So comparison of MRL and HOS with FE depends on the
solution, 
but in any case 
accuracy of HOS it contingent on $\mA_0$, $\mA_1$ and
$\mA_2$ not being large at the same time, to ensure relative smallness
of $\ErrOS$.

\subsection{Implementation}

Most of the algorithms described here were implemented in C language
in double precision floating point arithmetics and compiled using GNU
Compiler Collection (GCC) (version 4.7.2). The exception is
computation of eigenvalues and eigenvector tables, which was done
using mathematical software Sage~\cite{sage} (version 5.9).
Simulation were performed on Intel Core i5-3470 CPU with the clock
frequency 3.20GHz under GNU/Linux operating system (distribution
Fedora 18).

\section{Results}

\begin{figure*}
  \centering
  \includegraphics[width=7.16in]{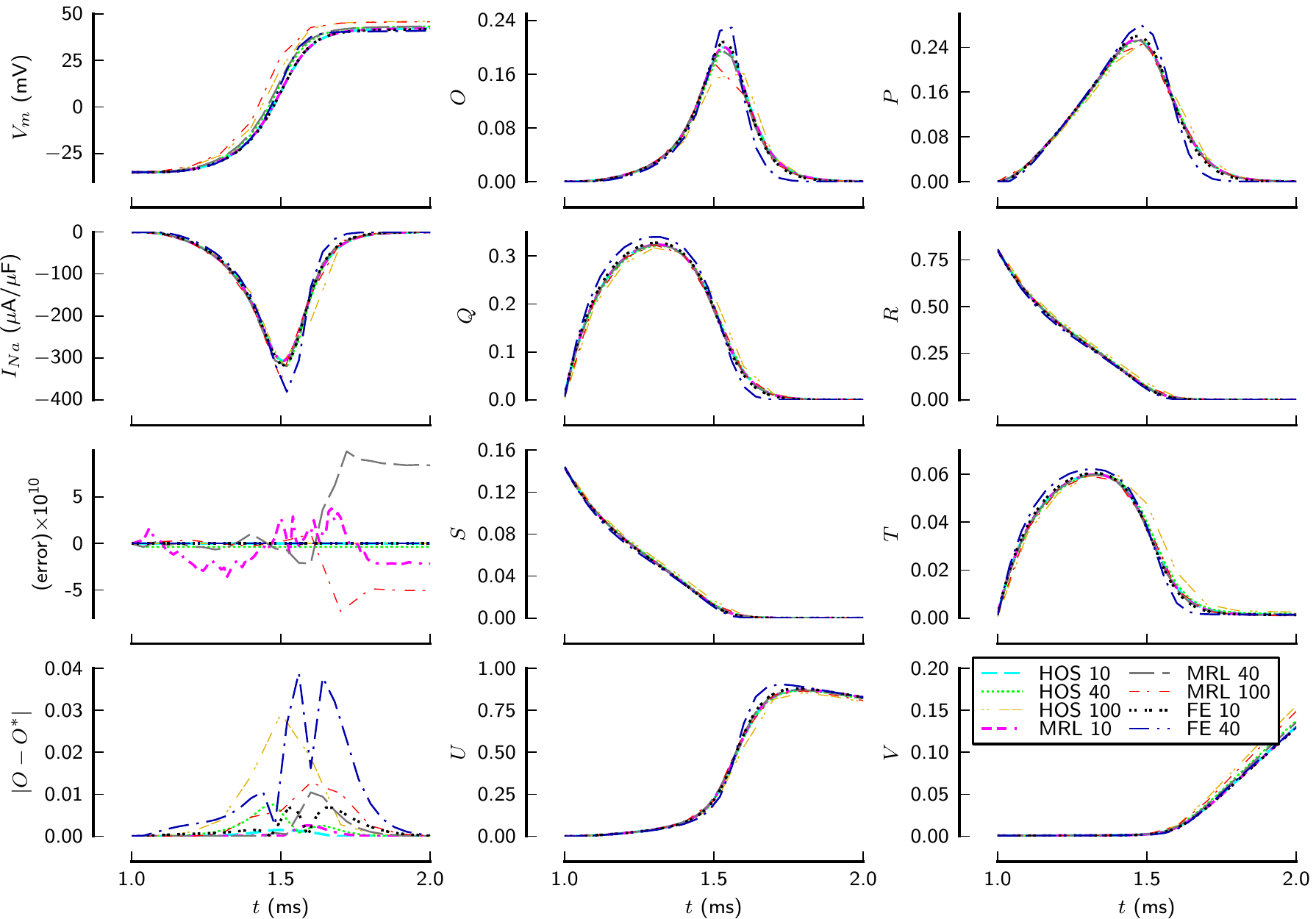}
  \caption[]{%
    Cardiac cell action potential (AP) simulations (detail of the
    first second) with $\INa$ Markov chain model solved using: %
    hybrid operator splitting method (HOS) -- light colours; %
    matrix Rush-Larsen (MRL) -- semi-dark colours; %
    and forward Euler (FE) -- dark colours; %
    The left column shows: %
    membrane potential ($\Vm$), sodium current ($\INa$), error
    calculated as deviation from the states conservation law, %
    and deviation of the state occupancy $\sO$ from the solution
      using FE with $\dt=1\,\us$ ($\sO^*$); %
    the other two columns show the state occupancies. %
    The model was solved with time step $\dt=10\,\us$, $\dt=40\,\us$,
    and $\dt=100\,\us$ represented by thick, middle thick, thin lines
    respectively.%
  }
  \figlabel{siml}
\end{figure*}

\Fig{siml} shows the detail of the first millisecond of simulated
cardiac excitation, the onset of an action potential (AP).  The $\INa$
Markov chain model was solved using the three suggested integration
methods: forward Euler (FE), matrix Rush-Larsen (MRL), and
hybrid operator splitting (HOS), as described in the methods section.
The model was solved with time step $\dt=10\,\us$, $\dt=40\,\us$, and
$\dt=100\,\us$, except for FE, which was also solved for
  $\dt=1\,\us$, to be used as a reference, but not for $\dt=100\,\us$,
  due to instability.

The model excitable cell was  initially at the resting
state, and at the time $\t=1\,\ms$, an AP was initialized by %
instantaneous injection of potassium ions, raising
  the membrane potential %
to $\Vm=-35\,\mV$. The initial conditions of
the states of the $\INa$ MC model are specified in \tab{states-term}
and the initial states of the remaining variables of the model can be
found in the supplementary material. Before the initiation, more than
90\% of the channels reside in the states $\sR$ and $\sS$, which
require at least three transitions to get to the open state. After
the initiation, the channels start to transit rapidly
towards the open state $\sO$ and then to the state $\sU$. Within about
$0.7\,\ms$ almost all channels reside in the state $\sU$. Then, the
channels slowly transit to the state $\sV$, where they stay until the
resting potential $\Vm$ is recovered. The states $\sU$ and $\sV$ are
similar to the situation when the inactivation gate $\hgate$ is closed
in the gate model. The states $\sS$ and $\sT$ have less than 10\%
occupancy during all the stages of the action potential.%
The plot of $\sW$ is omitted, as this variable
  changes very little during the time interval shown.

The results for the time step $\dt=10\,\us$ are consistent in all
panels. The FE is still stable at time step
$\dt=40\,\us$, however, compared to MRL and HOS, the
FE solution is less accurate, resulting in a higher peak and faster
decay of both the open state $\sO$ occupancy and the resulting $\INa$
current.

Comparison of the solutions for $\sO(\t)$ with the
  reference $\Oref(\t)$, obtained by FE with $\dt=1\,\us$, is shown on
  \fig{siml} (first column, fourth row). %
We see that MRL and HOS
    approximate $\Oref$ better than FE at the same time steps.  This
    is consistent with results of evaluation of the error estimates
    over the AP solution: we have $\max(\ErrFE)\approx2700$,
    $\max(\ErrMRL)\approx118$, $\max(\ErrHOS)\approx125$, with
    $\max(\ErrOS)\approx19$, all in $\ms^{-2}$, and
    $\min(\ErrFE/\ErrHOS)\approx2.3$, $\min(\ErrFE/\ErrMRL)\approx3.2$
    (see the Supplement). %
This suggests that exponential integrators can be useful,
  for their accuracy, even when instability is not a concern, say in
  systems with slower dynamics, such as $\IKs$.

At longer time steps, FE is unstable
(\fig{stiffness} illustrates a mild case of the instability), while
MRL and HOS continue to provide stable solutions. At
$\dt=100\,\us$,
the peak of the most important component of 
$\vu{}$, the occupancy of the open state $\sO$, is slightly lower than 
at shorter timesteps. On the other hand, the decrease of the peak of
the total $\INa$ current in these two methods (HOS, MRL at $\dt =
100\,\us$) is relatively small compared to the decrease of the open
state occupancy. Also, the decay of the $\INa$ current in the MRL $\dt
= 100\,\us$ is slower than in the other cases.  Note that the lead of
the APs onsets at $\dt=100\,\us$ against smaller time steps is
comparable to the value of $\dt$. 

Approximation of the whole APs rather than just their onsets
  is illustrated in \fig{stable}.

Further increase of the time steps (not shown) in MRL and
  HOS gives significant errors in the AP, e.g. at $\dt=200\,\us$ there
  is a $30\,\mV$ overshoot. Stability persists for much longer: for
  HOS the solution becomes unphysical (a negative concentration) at
  about $\dt=2\,\ms$ without loss of stability, and for MRL an
  instability occurs at about $7.5\,\ms$, although the solution is
  then also very different from the true AP.

\begin{table}[t!]
  \caption[]{Elapsed simulation time [s]. Cell model \cite{Clancy-Rudy-2002} 100 pulses with CL=1000\,\ms.}
  \centering
    \begin{tabular}{l||cc|cc|cc}
      \hline 
      &\multicolumn{2}{c|}{$\dt=10\,\us$} &\multicolumn{2}{c|}{$\dt=40\,\us$} &\multicolumn{2}{c}{$\dt=100\,\us$} \\
      \cline{2-7}
     $\INa$ Model  &$\INa$ & Total &$\INa$ & Total &$\INa$ & Total \\
      \hline
      FE          & 4.88 & 22.34 & 1.24 & 5.59  &      &  \\
      FE (tab.)   & 2.48 & 19.98 & 0.60 & 5.01  &      &  \\      
      MRL (tab.)  & 2.96 & 20.45 & 0.74 & 5.16  & 0.28 & 2.06 \\
      HOS        & 8.11 & 25.71 & 2.01 & 6.43  & 0.81 & 2.58 \\
      HOS (tab.) & 2.81 & 20.31 & 0.71 & 5.11  & 0.29 & 2.05 \\
      \hline
    \end{tabular} \label{tab:cost}
\end{table}

\Tab{cost} illustrates the efficiency of the three methods at three
different time steps $\dt$. This was done by measuring time taken by
simulations consisting of 100 pulses with a cycle length (CL) of
$1000\,\ms$ without any output. The pulses were initialized by %
  an instantaneous injection of potassium ions of a sufficient amount to set
the membrane potential to $\Vm=-35\,\mV$. The table shows times taken
by the whole cell model (``Total'') and by the Markov Chain model
computations (``$\INa$'').  The times shown are median values from six
separate simulations in each case to minimize the effect of other
processes running on the computer.

At $\dt=10\,\us$, FE is the most efficient method. 
Computation of $\INa$ 
accounts for 21.8\% and 31.5\% of the overall
computation cost in FE and HOS respectively. Tabulation allows
reduction of the computation cost of $\INa$ by 49.1\% and 65.3\% in FE
and HOS. MRL was used only with tabulation using the precomputed
eigenvalues and eigenvectors matrices and the computational cost at
the $\dt=10\,\us$ is comparable with FE. These proportions are
consistent with the results at time step $\dt = 40\,\us$ and $\dt =
100\,\us$ for MRL and HOS.  So, at the same time step, the
computational costs of the proposed methods are slightly higher, but
the accuracies are somewhat better, compared to FE. The most important
benefit of HOS and MRL is, however, the possibility of using larger
time steps.

\section{Conclusion}

Both proposed methods maintain stability at larger time steps, and
improve the accuracy of the solution at the same time
  step, compared to the explicit ODE solver (FE). When tabulated,
those methods are comparable to FE in computational cost. As expected,
using larger time steps results in reduction of computational cost.

\begin{figure}
  \centering
  \includegraphics[width=3.5in]{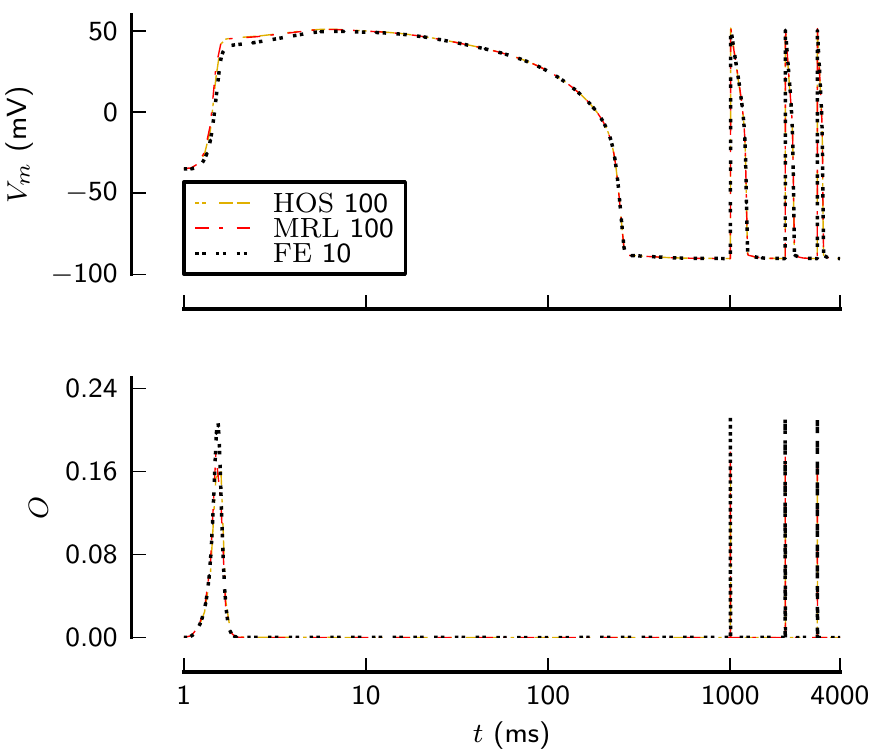}
  \caption[]{%
    Approximation properties of HOS ($\dt = 100\,\us$) and
    MRL ($\dt = 100\,\us$) compared to FE ($\dt = 10\,\us$)%
     during a single cell simulation of 4 action potentials with cycle length (CL) of 1000~ms (logarithmic scale): %
    membrane potential $\Vm$ (top panel), %
    open state occupancy $\sO$ (bottom panel).%
  }
  \figlabel{stable}
\end{figure}

MRL method extends the popular RL method, developed for gate models, to
Markov chain models. MRL is more universal than HOS, and
may be made ``automatic''. The only restriction of our implementation
is the assumption of diagonalizability of matrix $\mA(\Vm)$ for all
voltages.  If in another model this happens not to be the case, then
some more sophisticated approach would be needed. If non-diagonalizability is a
regular feature, say due to identical definitions of some of the TRs,
then a Jordan form can be used instead; if it only happens
at selected voltages, then interpolation of matrices $\mT{}(\Vm)$ may
be sufficient.

HOS method depends on the
possibility to split the transition rates to multiple (three in 
our case) sub-systems according to their speeds,
and solve each of the subsystems on its own. Our solution
benefits from the possibility of solving the fast subsystems
analytically. Implementing the analytical solution results in even
better speed-up as the resulting timestepping matrices are
sparse. However, the possibility of a suitable analytical
solution is not guaranteed for a general MC model. In this case, the
fast time subsystems can be solved using diagonalization like in the
MRL method, which might require additional computational time.

Finally we comment on the order of approximation. In this paper
  we considered first order schemes, and they are most popular in
  practice.  However, the approximation order can be improved by
using more sophisticated methods, both for the whole cell model (say
using Runge-Kutta approach) and for the exponential solvers. For the
original Rush-Larsen scheme, higher-order variants have been
  proposed and tested~\cite{Sundnes-etal-2009,Perego-Veneziani-2009},
  and the same ideas can be extended to the matrix case as well.
Naturally, HOS method may then need to involve a more
sophisticated operator splitting method to correspond.

Another appealing direction for further research is application of the
proposed methods to other important MC models. MRL is straightforward for
any MC where TRs depend only on one variable, otherwise tabulation
will be a bit more problematic. The success of the HOS
approach will depend on the asymptotic properties of the TRs.

\section*{Acknowledgment}

During this study, TS was funded by University of Exeter PhD
Studentship and VNB was partly supported by EPSRC grant
EP/I029664/1. We are grateful to C.E.~Clancy for the permission to use
the original authors code for this study and S.~Sherwin for
encouraging discussions.

\ifCLASSOPTIONcaptionsoff
  \newpage
\fi

\bibliographystyle{IEEEtran}

\cleardoublepage
\markboth{}{}

\textbf{\LARGE\centering
Supplementary material: \\ Exponential integrators for a Markov chain model \\ of the fast
  sodium channel of cardiocytes\\}
\vspace*{1.5\baselineskip}
\centerline{%
      Tom\' a\v s Star\' y, %
    Vadim N. Biktashev %
}
\vspace*{3.5\baselineskip}

The numeration of equations in this document continues from the main
text, and the literature references are to the literature list in the
main text, repeated in the end of this document for the reader's
convenience.

\section{Cell Model Definition}
%

This section contains the definition of the model according to the
authors' code~\cite{Clancy-Rudy-2002}. The format of equations and
subsections aims to correspond to the papers where those equations
were published to facilitate a straightforward comparison. The known
differences with the papers are marked by the sign:$\diff$.  
Voltages are measured in $\mV$, time in $\ms$ and concentrations in
mmol/L.%
The membrane currents are adjusted to the specific
  membrane capacitance $\Cm = 1\,\uF/\cm^2$\cite{Luo-Rudy-1994a} and
  are measured in $\uA/\uF$.

\subsection*{Standard ionic concentrations}
\begin{align}
\MACNaout =& 140\eqlabel{MACNaout}\diff  \\ %
\MACKout  =& 4.5 \eqlabel{MACKout}  \diff  \\ %
\MACcaout =& 1.8
\end{align}
which differs from~\cite{Luo-Rudy-1994a} where $\MACNaout = 150$;
$\MACKout = 5.4$.

\subsection*{Initial Values of Variables}
\begin{align}
\MACxsOne =& 0  \diff \\
\MACxsTwo =& 0  \diff \\
\MACv =& -95 \\ 
\MACcansr =&	1.8 \\ %
\MACcajsr =&	1.8 \\ %
\MACcamyo =& 0.00012 \\ %
\MACb =& 0.00141379 \\ %
\MACg =& 0.98831 \\ %
\MACd   =& 6.17507\D{-6} \\ %
\MACf   =& 0.999357 \\ %
\MACxr	       =& 2.14606\D{-4} \\ %
\MACNain  =&	7.9 \diff \\ %
\MACKin   =& 147.23 \eqlabel{MACKin}\diff 
\end{align}
which differs from~\cite{Luo-Rudy-1994a} where $\MACKin = 145$;
$\MACNain  =	10$, and no initial values were given for
$\MACxsOne$ and 
$\MACxsTwo$.

\subsection*{Physical Constants}
\begin{align}
  \MACr =& 8314 \\ %
\MACF     =&	96485 \\ %
\MACT     =& 310
\end{align}

\subsection*{Cell geometry}
\begin{align}
\MACl     =& 0.01 \\ %
\MACrad   =& 0.0011 \\ %
\MACvcell =& 3.801\D{-5} \\ %
\MACageo  =& 2 \MACpi \MACrad ^2+2 \MACpi \MACrad \MACl \\ %
\MACacap  =& 2\MACageo  \\ %
\MACvmyo  =& 2.58468\D{-5} \\ %
\MACvnsr  =& 0.0552\MACvcell \\ %
\MACvjsr  =& 0.0048\MACvcell
\end{align}

\subsection*{$\Naion^+\!$-$\Kion^+$ pump : $\MACInak$}
\begin{align}
  \MACInak =& 1.5 \MACfNaK \frac{1}{1+(10/\MACNain)^{1.5}} \cdot  \frac{\MACKout}{\MACKout+1.5} \\ %
  \MACfNaK =& \frac{1}{1+0.1245 \MACexp\left(-0.1\cdot \frac{\MACv \MACF}{\MACr \MACT}\right)+0.0365 \MACsig \MACexp((-\MACv \MACF)/(\MACr \MACT))} \\ %
  \MACsig =& \frac{1}{7}\MACexp\left(\frac{\MACNaout}{67.3}\right)-1
\end{align}
which is identical to~\cite{Luo-Rudy-1994a}

\subsection*{$\MACiks$, the Slow Component of the Delayed Rectifier $\Kion^+$ Current}
\begin{align}
  \MACiks =& \MACGks \MACxsOne \MACxsTwo (\MACv-\MACEks) \\ %
  \MACEks =& (\MACr \MACT/\MACF) \MAClog((4.5+\MACprnak 150)/(\MACKin+\MACprnak \MACNaout)) \eqlabel{MACEks-code}\diff \\ %
  \MACprnak =& 0.01833 \\ %
  \MACGks =& (0.433 (1+0.6/(1+(0.000038/\MACcamyo)^{1.4})))\cdot 0.615 \eqlabel{MACGks-code}\diff 
\\ 
%
  \MACxsOness	 =& 1/(1+\MACexp(-(\MACv-1.5)/16.7)) \\ %
  \MACxsTwoss	 =& \MACxsOness \\ %
  \MACtauxsOne =& \left(0.0000719 \frac{\MACv+30}{1-\MACexp(-0.148 (\MACv+30))}+0.000131 \frac{\MACv+30}{\MACexp(0.0687 (\MACv+30))-1}\right)^{-1} \\ 
  \MACtauxsTwo =& 4 \MACtauxsOne
\end{align}
The definition of $\MACEks$ in equation \eq{MACEks-code} differs 
from~\cite{Viswanathan-etal-1999a} by the hard-coded
term for the $\MACKout=4.5$ 
and $\MACNaout=150$ rather
than values defined by equations \eqtwo(MACNaout,MACKout)
where $\MACKout=4.5$ and $\MACNaout=140$ are parameters.

The definition of $\MACGks$ in equation \eq{MACGks-code} is
multiplied by 0.615 ``to simulate the intramural heterogeneity'', %
which is slightly different from the factor $0.652$ used in
  Viswanathan et al. (1999)~\cite{Viswanathan-etal-1999a} to simulate
  epicardial cell.
 
Otherwise the $\MACiks$ definition is identical
to~\cite{Viswanathan-etal-1999a}.

\begin{align}
  \ddt{\MACxsOne} =& \frac{\MACxsOness - \MACxsOne}{\MACtauxsOne} \\ %
  \ddt{\MACxsTwo} =& \frac{\MACxsTwoss - \MACxsTwo}{\MACtauxsTwo} 
\end{align}

\subsection*{$\MACikr$, the Fast Component of the Delayed Rectifier
  $\Kion^+$ Current}
\begin{align}
  \MACikr =& \MACgkr \MACxr \MACrot (\MACv-\MACekr) \\ %
  \MACgkr =& 0.02614 \sqrt{\MACko/5.4} \\ %
  \MACxrss	=& 1/(1+\MACexp(-(\MACv+21.5)/7.5)) \\ %
  \MACrot =& 1/(1+\MACexp((\MACv+9)/22.4)) \\ %
  \MACekr =& ((\MACR \MACtemp)/\MACfrdy) \MAClog(\MACko/\MACki) \\ %
  \MACtauxr =& \left(0.00138 \frac{\MACv+14.2}{1-\MACexp(-0.123 (\MACv+14.2))}+0.00061 \frac{\MACv+38.9}{\MACexp(0.145 (\MACv+38.9))-1}\right)^{-1}
\end{align}
which is identical to~\cite{Zeng-etal-1995}. The original notation for
$\MACrot$ was $\sR$; we use $\MACr$ for the gas constant in this
  section, and for one of the state occupancies of the Markov Chain
  model elsewhere in the rest of the paper.

\begin{align}
  \ddt{\MACxr} =& \frac{\MACxrss-\MACxr}{\MACtauxr}
\end{align}

\subsection*{ Time-independent $\Kion^+$ current: $\MACIkOne$}
\begin{align}
  \MACIkOne =& \MACGkOne \MACkOne (\MACv-\MACEkOne) \\ %
  \MACEkOne =& (\MACr \MACT/\MACF) \MAClog(\MACKout/\MACKin) \\ %
  \MACGkOne =&	0.75\cdot \sqrt{(\MACKout/5.4)} \\ 
  \MACakOne =& 1.02/(1+\MACexp(0.2385 (\MACv-\MACEkOne-59.215))) \\ %
  \MACbkOne =&	\frac{0.49124 \MACexp(0.08032 (\MACv-\MACEkOne+5.476))+\MACexp(0.06175 (\MACv-\MACEkOne-594.31))}{1+\MACexp(-0.5143 (\MACv-\MACEkOne+4.753))}
\end{align}
which is identical to~\cite{Luo-Rudy-1994a}.

\begin{align}
  \MACkOne  =& \MACakOne/(\MACakOne+\MACbkOne) \\ %
\end{align}

\subsection*{Plateau $\Kion^+$ current: $\MACIkp$}
\begin{align}
  \MACIkp =& 0.00552 \MACkp (\MACv-\MACEkOne) \\ %
  \MACkp  =&	1/(1+\MACexp((7.488-\MACv)/5.98)) 
\end{align}
equivalent to~\cite{Luo-Rudy-1994a} with an update from~\cite{Zeng-etal-1995}.
\begin{align}
  \MACik =& \MACIkOne+\MACIkp
\end{align}

\subsection*{Currents through the L-type $\Caion^{+2}$ channel $\MACicaL$}
\begin{align}
  \MACicaL =& \MACilca+\MACilcak+\MACilcana \\ %
  \MACilca	 =& \MACd \MACf \MACfca \MACibarca \\ %
  \MACilcak	 =& \MACd \MACf \MACfca \MACibark \\ %
  \MACilcana =& \MACd \MACf \MACfca \MACibarna \\ %
  \MACibarca =& \MACpca  \MACzca^2 \frac{(\MACv \MACfrdy^2)}{\MACR \MACtemp}\,\frac{\MACgacai \MACcai \,\MACexp((\MACzca \MACv \MACfrdy)/(\MACR \MACtemp))-\MACgacao \MACcao}{\MACexp((\MACzca \MACv \MACfrdy)/(\MACR \MACtemp))-1} \\ %
  \MACibarna =& \MACpna \MACzna^2 \frac{(\MACv \MACfrdy^2)}{\MACR \MACtemp}\,\frac{\MACganai \MACnai \,\MACexp((\MACzna \MACv \MACfrdy)/(\MACR \MACtemp))-\MACganao \MACnao}{\MACexp((\MACzna \MACv \MACfrdy)/(\MACR \MACtemp))-1} \\ %
  \MACibark =& \MACpk \MACzk^2 \frac{(\MACv \MACfrdy^2)}{\MACR \MACtemp}\,\frac{\MACgaki \MACki \,\MACexp((\MACzk \MACv \MACfrdy)/(\MACR \MACtemp))-\MACgako \MACko}{\MACexp((\MACzk \MACv \MACfrdy)/(\MACR \MACtemp))-1} \\ %
  \MACpca   =& 5.4\D{-4} \hspace{2em}
  \MACgacai = 1 \hspace{2em}
  \MACgacao = 0.341 \\ %
  \MACpna   =& 6.75\D{-7}  \hspace{2em}
  \MACganai =	0.75 \hspace{2em}
  \MACganao =	0.75 \\ %
  \MACpk    =&	1.93\D{-7}  \hspace{2em}
  \MACgaki  =	0.75 \hspace{2em}
  \MACgako  =	0.75 \\ %
  \MACfca =& 1/(1+\MACcai/\MACkmca) \\ %
  \MACkmca  =&	0.0006 \\ %
  \MACdss  =& 1/(1+\MACexp(-(\MACv+10)/6.24)) \\ %
  \MACtaud =& \MACdss (1-\MACexp(-(\MACv+10)/6.24))/(0.035 (\MACv+10)) \\ %
  \MACfss  =& (1/(1+\MACexp((\MACv+32)/8)))+(0.6/(1+\MACexp((50-\MACv)/20))) \eqlabel{MACfss-code}\diff \\ %
  \MACtauf =& 1/(0.0197 \MACexp(-(0.0337 (\MACv+10)^2))+0.02)
\end{align}
Equation~\eq{MACfss-code} differs from~\cite{Luo-Rudy-1994a}
which has $8.6$ rather than $8$ in the denominator of the argument of
the first exponential. Otherwise, these equations are exactly the same
as in~\cite{Luo-Rudy-1994a}.
\begin{align}
  \MACzna =& 1 \\ %
  \MACzk  =& 1 \\ %
  \MACzca =& 2 \\ %
  \ddt{\MACd} =& \frac{\MACdss-\MACd}{\MACtaud} \\ %
  \ddt{\MACf} =& \frac{\MACfss-\MACf}{\MACtauf}
\end{align}

\subsection*{$\Caion^{2+}$ Current Through \sT-Type $\Caion^{2+}$ Channels $\MACicaT$}
\begin{align}
  \MACicaT =& \MACgcat  \MACb^2 \MACg (\MACv-\MACeca) \\ %
  \MACgcat =& 0.05 \\ %
  \MACbss =& 1/(1+\MACexp(-(\MACv+14)/10.8)) \\ %
  \MACggs =& 1/(1+\MACexp((\MACv+60)/5.6)) \\ %
  \MACeca =& (\MACr \MACT/(2 \MACF)) \, \MAClog(\MACcaout/\MACcamyo) \\ %
  \MACtaub =& 3.7+6.1/(1+\MACexp((\MACv+25)/4.5)) \\ %
    \MACtaug =& -0.875 \MACv+12 \text{ for: } \MACv \leq 0 \text{; and }
     \MACtaug = 12  \text{ for: } \MACv > 0
\end{align}
which correspond exactly to~\cite{Zeng-etal-1995}.
\begin{align}
  \ddt{\MACb} =& \frac{\MACbss-\MACb}{\MACtaub} \\ %
  \ddt{\MACg} =& \frac{\MACggs-\MACg}{\MACtaug}
\end{align}

\subsection*{$\Naion^+$-$\Caion^+$ exchanger: $\MACInaca$}
\begin{align}
  \MACInaca =&  \frac{2.5\D{-4}\MACexp((\MACgammas-1) \MACv \frac{\MACF}{\MACR \MACT}) (\MACexp(\MACv \frac{\MACF}{\MACR \MACT})  \MACnai^3 \MACcao-\MACnao^3 \MACcai)}{1+1\D{-4} \MACexp((\MACgammas-1) \MACv \frac{\MACF}{\MACR \MACT}) (\MACexp(\MACv \frac{\MACF}{\MACR \MACT}) \MACnai^3 \MACcao+\MACnao^3 \MACcai)} \diff \\ %
  \MACgammas =& 0.15 \diff
\end{align}
Here $\MACInaca$ depends on external $\MACcao,\MACnao$ as well as
internal $\MACnai,\MACcai$ concentrations, which is different
from~\cite{Luo-Rudy-1994a} where it depended only on external concentrations
$\MACcao,\MACnao$. The variable $\MACgammas = 0.35$ in~\cite{Luo-Rudy-1994a}.

\subsection*{Nonspecific $\Caion^{2+}$-activated current: $\MACinsca$}
\begin{align}
  \MACInsk =& 1.75\D{-7} \frac{\MACv  \MACF^2}{\MACr \MACT} \cdot \frac{0.75 \MACKin \MACexp(( \MACv \MACF)/(\MACr \MACT))-0.75 \MACKout}{\MACexp( \MACv \MACF/(\MACr \MACT))-1} \\ %
  \MACikns =& \MACInsk \frac{1}{1+(0.0012/\MACcamyo)^3} \\ %
  \MACInsna	=& 1.75\D{-7} \frac{\MACv \MACF^2}{\MACr \MACT}\cdot \frac{0.75 \MACNain \MACexp(( \MACv \MACF)/(\MACr \MACT))-0.75 \MACNaout}{\MACexp( \MACv \MACF/(\MACr \MACT))-1} \\ %
  \MACinsna	=& \MACInsna \frac{1}{1+(0.0012/\MACcamyo)^3 } \\ %
  \MACinsca	=& \MACikns+\MACinsna \\ %
  \MACPnsca	=& 1.75\D{-7}
\end{align}
This is almost identical to~\cite{Luo-Rudy-1994a} except the latter also made a
definition for $E_{ns(\Caion)}$ which however was not used.

\subsection*{Sarcolemmal $\Caion^{+2}$ pump: $\MACipca$}
\begin{align}
  \MACipca =& 1.15 \frac{\MACcamyo}{0.0005+\MACcamyo}
\end{align}
identical to~\cite{Luo-Rudy-1994a}.

\subsection*{$\Caion^{+2}$ background current: $\MACicab$}
\begin{align}
  \MACicab =& 0.003016 (\MACv-\MACEca) \\ %
  \MACEca  =& \MACr \MACT/(2 \MACF) \MAClog(\MACcaout/\MACcamyo)
\end{align}
identical to~\cite{Luo-Rudy-1994a}.

\subsection*{$\Naion^+$ background current: $\MACinab$}
\begin{align}
  \MACEna =& ((\MACR \MACtemp)/\MACfrdy)\log(\MACnao/\MACnai) \\
  \MACinab =& 0.00141 (\MACv-\MACEna)
\end{align}
identical to~\cite{Luo-Rudy-1994a}.

\subsection*{$\Caion^{2+}$ uptake and leakage of NSR: $\MACIup$ and $\MACIleak$}
\begin{align}
  \MACIup	=& 0.00875 \MACcamyo/(\MACcamyo+0.00092) \diff \\ %
  \MACKleak	=& 0.005/15 \\ %
  \MACIleak	=& \MACKleak \MACcansr
\end{align}
The definition of $\MACIup$ in~\cite{Luo-Rudy-1994a} is
ambiguous. This version is consistent with one possible understanding.

\subsection*{ $\Caion^{+2}$ Fluxes in NSR}

\begin{align}
  \ddt{\MACcansr} =&  (\MACIup-\MACIleak-\MACItr \MACvjsr/\MACvnsr)
\end{align}

\subsection*{$\Caion^{2+}$ Fluxes in Myoplasm}
  \begin{align}
  \MACitca =& \MACilca+\MACicab+\MACipca-2 \MACInaca+\MACicaT \\ %
  \MACdcai =& -\MACdt (((\MACitca \MACacap)/(\MACvmyo 2 \MACF))+((\MACIup-\MACIleak) \MACvnsr/\MACvmyo)-(\MACIrel \MACvjsr/\MACvmyo))
           \eqlabel{DeltaCai} \\ %
 \MACcaion	=& \MACTRPN+\MACCMDN+\MACdcai+\MACcamyo \\ %
  \MACB =& 0.05+0.07-\MACcaion+0.0005+0.00238 \\ %
  \MACC =& (0.00238\cdot 0.0005)-(\MACcaion (0.0005+0.00238))+(0.07\cdot 0.00238)+(0.05\cdot 0.0005) \\ %
  \MACD =& -0.0005\cdot 0.00238 \MACcaion \\ %
  \MACFab =& \sqrt{( \MACB^2-3 \MACC)} \\ %
  \MACcamyo	=& 1.5 \MACFab \MACcos(\MACacos((9 \MACB \MACC-2 \MACB^3-27 \MACD)/(2 (\MACB^2-3 \MACC)^{1.5}))/3)-(\MACB/3)
  \eqlabel{cubic-formula}
\end{align}
This definition merely summarises computations that are done in the
code, which \textit{de facto} describe a time-stepping algorithm for a
system of a differential equation and a finite constraint, rather than
the equation and the constraint themselves, hence the time step $\MACdt$ is
present in~\eq{DeltaCai}. Any attempts of higher-order numerical
approximations would have to take this into account.  The explicit
solution of the finite constraint given by the cubic
formula~\eq{cubic-formula} follows~\cite{Zeng-etal-1995}
whereas~\cite{Luo-Rudy-1994a} used Steffensen's iterations for that purpose.

\subsection*{$\Caion^{2+}$ Fluxes in JSR}
  \begin{align}
\MACdcajsr =& \MACdt (\MACItr-\MACIrel) \eqlabel{DeltaCaJSR} \\ %
 \MACbjsr	 =& 10-\MACCSQN-\MACdcajsr-\MACcajsr+0.8 \\ %
 \MACcjsr	 =& 0.8 (\MACCSQN+\MACdcajsr+\MACcajsr) \\ %
  \MACcajsr =& (\sqrt{( \MACbjsr^2+4 \MACcjsr)}-\MACbjsr)/2 
  \end{align}
Ditto: $\MACdt$ is present in~\eq{DeltaCaJSR}. 

\subsection*{Sodium Ion Fluxes}
  \begin{align}
  \MACitna =& \MACina+\MACinab+\MACilcana+\MACinsna+3 \MACInak+3 \MACInaca \\ %
  \ddt{\MACNain}	=& - (\MACitna \MACacap)/(\MACvmyo \MACF)
\end{align}
$\MACNain$ is constant in \cite{Luo-Rudy-1994a,Zeng-etal-1995,Viswanathan-etal-1999a}.

\subsection*{Potassium Ion Fluxes}

\begin{align}
  \MACitk =& \MACikr+\MACiks+\MACik+\MACilcak+\MACikns-2 \MACInak+\MACito+\MACstimulus \\ %
  \ddt{\MACKin} =& - (\MACitk \MACacap)/(\MACvmyo \MACF)
\end{align}
$\MACKin$  is constant in \cite{Luo-Rudy-1994a,Zeng-etal-1995,Viswanathan-etal-1999a}.

\subsection*{$\Caion ^{2+}$ buffers in the myoplasm}
\begin{align}
  \MACTRPN =& 0.07 \MACcamyo/(\MACcamyo+0.0005) \\ %
  \MACCMDN =& 0.05 \MACcamyo/(\MACcamyo+0.00238)
\end{align}
identical to~\cite{Luo-Rudy-1994a}.

\subsection*{$\Caion ^{2+}$ buffer in JSR and SCQN}
\begin{align}
  \MACCSQN =& 10 (\MACcajsr/(\MACcajsr+0.8)) \\ %
\end{align}
identical to~\cite{Luo-Rudy-1994a}.

\subsection*{CICR From Junctional SR (JSR)}

\begin{align}
  \MACIrel =& \MACGrel \MACryropen \MACryrclose (\MACcajsr-\MACcamyo) \\ %
  \MACGrel =& 150/(1+\MACexp(\MACitca+5)/0.9) 
                                \diff \eqlabel{Grel} \\ %
  \MACryropen =& 1/(1+\MACexp((-\MACtc+4)/0.5)) 
                                \diff \eqlabel{ryropen} \\ %
  \MACryrclose =& 1-(1/(1+\MACexp((-\MACtc+4)/0.5))) 
                                \diff \eqlabel{ryrclose}
\end{align}
Here is another deviation of the model description from the standard
form of a system of ODEs, and this also would have to be taken into
account in any attempts of higher-order schemes.  Variables
$\MACryrclose$ and $\MACryropen = 1-\MACryrclose$ ensure that the
calcium release channel is open at a fuzzy time interval around 4~ms
after the steepest point of the upstroke of the action potential. This is
done using an additional time variable $\MACtc$ which is linked to the
$\MACt$, that is $\d{\MACtc}/\d{\t}=1$ most of the time, except $\MACtc$ is
reset to zero each time the $\MACdvdt$ reaches a \emph{significant}
local maximum, ``significant'' meaning $\MACdvdt>1\,\mV/\ms$. 
In~\cite{Luo-Rudy-1994a}, $\MACGrel$ is defined differently from \eq{Grel}
and calcium release proceeds with a different dynamics
from~\eqtwo(ryropen,ryrclose), e.g. it starts sharply $2\,\ms$ after
the the time of the maximum $\MACdvdt$.

\subsection*{Translocation of $\Caion^{2+}$ ions from NSR to JSR:
  $\MACItr$ }
\begin{align}
  \MACItr =& (\MACcansr-\MACcajsr)/180
\end{align}
identical to~\cite{Luo-Rudy-1994a}. 

\subsection*{Total time-independent current: $\MACiv$}
\begin{align}
  \MACiv =& \MACinab+\MACInak+\MACipca+\MACIkp+\MACicab+\MACIkOne
\end{align}
identical to~\cite{Luo-Rudy-1994a}.

\subsection*{Total Current}
\begin{align}
  \MACit =& \MACikr+\MACiks+\MACik+\MACilcak+\MACikns-2 \MACInak+\MACina+\MACinab+\MACilcana+\MACinsna+3 \MACInak+ \nonumber \\&3 \MACInaca+\MACilca+\MACicab+\MACipca-2 \MACInaca+\MACicaT
  \end{align}

\subsection*{ Membrane Potential}
  \begin{align}
    \ddt{\MACv} =& -\MACit .
  \end{align}

\section{$\INa$ Markov Chain Model Definition}

Up to the choice of notation, we use the model described
in~\cite{Clancy-Rudy-2002}. The fast sodium current is defined by
\begin{align}
  \MACina =  \MACGNa (\MACv-\MACEna) \sO , \eqlabel{MACina-code}
\end{align}
where the channel open probability $\sO$ is defined by the system of ODEs
\begin{align}
 \ddt{\sO} =& \alPO \sP + \alUO \sU - (\alOP+\alOU) \sO                 \eqlabel{dINa-O-new} \\
 \ddt{\sP} =& \alQP \sQ + \alUP \sU + \alOP \sO - (\alPQ+\alPU+\alPO) \sP \eqlabel{dINa-P-new} \\
 \ddt{\sQ} =& \alRQ \sR + \alTQ \sT + \alPQ \sP - (\alQR+\alQT+\alQP) \sQ \eqlabel{dINa-Q-new} \\
 \ddt{\sR} =& \alSR \sS + \alQR \sQ - (\alRS+\alRQ) \sR                 \eqlabel{dINa-R-new} \\
 \ddt{\sS} =& \alTS \sT + \alRS \sR - (\alST+\alSR) \sS                 \eqlabel{dINa-S-new} \\
 \ddt{\sT} =& \alQT \sQ + \alST \sS + \alUT \sU - (\alTQ+\alTS+\alTU) \sT \eqlabel{dINa-T-new} \\
 \ddt{\sU} =& \alTU \sT + \alPU \sP + \alVU \sV + \alOU \sO - (\alUT+\alUP+\alUO+\alUV) \sU 
                                                                \eqlabel{dINa-U-new} \\
 \ddt{\sV} =& \alUV \sU + \alWV \sW - (\alVU+\alVW) \sV                 \eqlabel{dINa-V-new} \\
 \ddt{\sW} =& \alVW \sV - \alWV \sW                                   \eqlabel{dINa-W-new} 
\end{align}
with the transition rates defined by
\begin{align*}
\alRQ=\alST=\alpha_{11} =& \frac{3.802}{0.1027 \,\e^{-\Vm/17.0} + 0.20 \,\e^{-\Vm/150}}  \\
\alQP=\alTU=\alpha_{12} =& \frac{3.802}{0.1027 \,\e^{-\Vm/15.0} + 0.23 \,\e^{-\Vm/150}}  \\
\alPO=\alpha_{13} =& \frac{3.802}{0.1027 \,\e^{-\Vm/12.0} + 0.25 \,\e^{-\Vm/150}}  \\
\alQR=\alTS=\beta_{11} =& 0.1917 \,\e^{-\Vm/20.3}  \\
\alPQ=\alUT=\beta_{12} =& 0.20 \,\e^{-(\Vm-5)/20.3}  \\
\alOP=\beta_{13} =& 0.22 \,\e^{-(\Vm-10)/20.3}  \\
\alUP=\alTQ=\alSR=\alpha_{3} =& 3.7933\cdot 10^{-7} \,\e^{-\Vm/7.7}  \\
\alPU=\alQT=\alRS=\beta_{3} =& 8.4\cdot 10^{-3} + 2\cdot 10^{-5} \Vm  \\
\alOU=\alpha_{2} =& 9.178 \,\e^{\Vm/29.68}  \\
\alUO=\beta_{2} =& \frac{\alpha_{13} \alpha_{2} \alpha_{3}}{\beta_{13} \beta_{3}}  \\
\alUV=\alpha_{4} =& \alpha_{2}/100  \\
\alVU=\beta_{4} =& \alpha_{3}  \\
\alVW=\alpha_{5} =& \alpha_{2} /(9.5 \cdot 10^4 )  \\
\alWV=\beta_{5} =& \alpha_{3} /50 
\end{align*}

\section{ Details of the hybrid Operator Splitting method}

\subsection{Operator splitting}
In the hybrid method we use operator splitting. The system of
equations~\eq{dINa-O-new}--\eq{dINa-W-new} is considered as an ODE
\begin{align*}
  \ddt{\vu{}} = \mA(\Vm(\t)) \vu{} , \tag{\ref{eq:markov-generic}}
\end{align*}
for the vector-function $\vu{}=\left(\sO,\sP,\sQ,\sR,\sS,\sT,\sU,\sV,\sW\right)\T=\vu{}(\t)$, and
the transition matrix is split into the sum
\begin{align}
  \mA = \mA_0 + \mA_1 + \mA_2
\end{align}
of the matrix $\mA_0$ of transition rates that are fast at high values
of $\Vm$,
\begin{align} 
  \mA_0 =&
  \begin{bmatrix}
    -\alOU & \alPO & 0      & 0      & 0      & 0      & 0 & 0 & 0\\
    0      & -\alPO & \alQP & 0      & 0      & 0      & 0 & 0 & 0\\
    0      & 0      & -\alQP & \alRQ & 0      & 0      & 0 & 0 & 0\\
    0      & 0      & 0      & -\alRQ & 0      & 0      & 0 & 0 & 0\\
    0      & 0      & 0      & 0      & -\alST & 0      & 0 & 0 & 0\\
    0      & 0      & 0      & 0      & \alST & -\alTU & 0 & 0 & 0\\
    \alOU & 0      & 0      & 0      & 0      & \alTU & 0 & 0 & 0\\
    0      & 0      & 0      & 0      & 0      & 0      & 0 & 0 & 0\\
    0      & 0      & 0      & 0      & 0      & 0      & 0 & 0 & 0\\
  \end{bmatrix},                                \eqlabel{A0}
\end{align}
the matrix $\mA_1$ of transition rates that are fast at low values of $\Vm$, 
\begin{align}
  \mA_1 =&
  \begin{bmatrix}
    -\alOP&	0&	0&	0&	0&	0&	0&	0&	0 \\
    \alOP&	-\alPQ&	0&	0&	0&	0&	0&	0&	0 \\
    0&	\alPQ&	-\alQR&	0&	0&	0&	0&	0&	0 \\
    0&	0&	\alQR&	0&	0&	0&	0&	0&	0 \\
    0&	0&	0&	0&	0&	\alTS&	0&	0&	0 \\
    0&	0&	0&	0&	0&	-\alTS&	\alUT&	0&	0 \\
    0&	0&	0&	0&	0&	0&	-\alUT&	0&	0 \\
    0&	0&	0&	0&	0&	0&	0&	0&	0 \\
    0&	0&	0&	0&	0&	0&	0&	0&	0 \\
  \end{bmatrix},                                \eqlabel{A1}
\end{align}
and the matrix $\mA_2$ of uniformly slow transition rates, 
\begin{align}
  \mA_2 =&
  \begin{bmatrix}
    0 & 0      & 0      & 0      & 0      & 0      & \alUO                & 0              & 0 \\
    0 & -\alPU & 0      & 0      & 0      & 0      & \alUP                & 0              & 0 \\
    0 & 0      & -\alQT & 0      & 0      & \alTQ  & 0                    & 0              & 0 \\
    0 & 0      & 0      & -\alRS & \alSR  & 0      & 0                    & 0              & 0 \\
    0 & 0      & 0      & \alRS  & -\alSR & 0      & 0                    & 0              & 0 \\
    0 & 0      & \alQT  & 0      & 0      & -\alTQ & 0                    & 0              & 0 \\
    0 & \alPU  & 0      & 0      & 0      & 0      & -(\alUP+\alUO+\alUV) & \alVU          & 0 \\
    0 & 0      & 0      & 0      & 0      & 0      & \alUV                & -(\alVU+\alVW) & \alWV \\
    0 & 0      & 0      & 0      & 0      & 0      & 0                    & \alVW          & -\alWV \\
  \end{bmatrix}.   	    	     	      	  \eqlabel{A2}
\end{align}

\subsection{The substeps}

Every timestep is then done in three substeps, each using one of the
three matrices $\mA_\m$, $\m=0,1,2$: 
\begin{align*}
  \vu{\n+1/3} =&\, \exp(\dt \mA_0(\Vm(\t_\n))) \, \vu{\n} , \tag{\ref{eq:sol-subs-1}} \\
  \vu{\n+2/3} =&\, \exp(\dt \mA_1(\Vm(\t_\n))) \, \vu{\n+1/3} , \tag{\ref{eq:sol-subs-2}} \\
  \vu{\n+1} =&\, \vu{\n+2/3} +\dt \mA_2(\Vm(\t_\n)) \, \vu{\n+2/3} . \tag{\ref{eq:sol-subs-3}}
\end{align*}
Note that $\Vm$ in all cases is evaluated at $\t=\t_\n$, in which we
simply follow the original Rush-Larsen idea of ``freezing'' $\Vm$ for
the duration of the time step. The matrix exponentials in
\eq{sol-subs-1} and \eq{sol-subs-2} can be understood in terms of
matrix Taylor series~\cite{Reidman-Trivedi-1988}, or the product of
the matrix exponential by the corresponding vector $\vu{}$ can be
understood just as the solutions of an initial-value problem for the
corresponding system of ODEs with constant coefficients.  The mapping
\eq{sol-subs-1} is calculated by solving the following initial-value
problem, defined by the matrix $\mA_0$ \eq{A0},
\begin{align}
  &\ddt{\sO} = - \alOU \sO + \alPO \sP, && \sO(0)=\sO_\n,\eqlabel{A0-O} \\
  &\ddt{\sP} = - \alPO \sP + \alQP \sQ, && \sP(0)=\sP_\n,\eqlabel{A0-P} \\
  &\ddt{\sQ} = - \alQP \sQ + \alRQ \sR, && \sQ(0)=\sQ_\n,\eqlabel{A0-Q} \\
  &\ddt{\sR} = - \alRQ \sR, 		  && \sR(0)=\sR_\n,\eqlabel{A0-R} \\
  &\ddt{\sS} = - \alST \sS,               && \sS(0)=\sS_\n,\eqlabel{A0-S} \\
  &\ddt{\sT} = \alST \sS - \alTU \sT,   && \sT(0)=\sT_\n,\eqlabel{A0-T} \\
  &\ddt{\sU} = \alOU \sO + \alTU \sT,   && \sU(0)=\sU_\n,\eqlabel{A0-U} \\
  &\ddt{\sV} = 0,                           && \sV(0)=\sV_\n,\eqlabel{A0-V} \\
  &\ddt{\sW} = 0,                           && \sW(0)=\sW_\n,\eqlabel{A0-W} 
\end{align}
and then evaluating the result at $\t=\dt$ to give $\sO_{\n+1/3}$,
$\dots$, $\sW_{\n+1/3}$. 
We note now that equations \eq{A0-R} and \eq{A0-S} are
decoupled and we can solve them to get
\begin{align}
  & \sR(\t)=\sR_\n\,\e^{-\alRQ \t}, \eqlabel{A0-Rsol}\\
  & \sS(\t)=\sS_\n\,\e^{-\alST \t}. \eqlabel{A0-Ssol} 
\end{align}
We then substitute \eq{A0-Rsol} into \eq{A0-Q} to obtain a closed
initial-value problem for $\sQ(\t)$,
\begin{align}
  \ddt{\sQ} + \alQP \sQ = \sR_\n\alRQ \,\e^{-\alRQ \t}, \quad \sQ(0)=\sQ_\n,
\end{align}
the solution of which is
\begin{align}
  \sQ(\t) = \sQ_\n \,\e^{-\alQP \t} - \sR_\n \frac{\alRQ (\,\e^{-\alQP
      \t} - \,\e^{-\alRQ \t})}{\alQP - \alRQ} . \eqlabel{A0-Qsol} \\
\end{align}
Similarly, we substitute \eq{A0-Ssol} into \eq{A0-T} to obtain 
\begin{align}
  \sT(\t) = \sT_\n \,\e^{-\alTU \t} - \sS_\n\frac{\alST(\,\e^{-\alTU\t} - \,\e^{-\alST\t})}{\alTU- \alST}. \eqlabel{A0-Tsol}\\
\end{align}
We then proceed in the same manner, by
substituting the obtained solution \eq{A0-Qsol} for
$\sQ(\t)$ into \eq{A0-P} to obtain $\sP(\t)$, and the solution
\eq{A0-Tsol} for $\sT(\t)$ into \eq{A0-U} to obtain $\sU(\t)$, and finally
the found solution for $\sP(\t)$ into \eq{A0-O} to obtain $\sO(\t)$. With
the obvious solutions to \eq{A0-V} and \eq{A0-W}, the result of all
these steps is mapping
\begin{align}
  &\sO_{\n+1/3} = \mOU \sO_\n + \kPO \sP_\n     + \kQO \sQ_\n  + \kRO \sR_\n ,\\
  &\sP_{\n+1/3} = \mPO \sP_\n + \kQP \sQ_\n + \kRP \sR_\n ,\\
  &\sQ_{\n+1/3} = \mQP \sQ_\n + \kRQ \sR_\n ,\\
  &\sR_{\n+1/3} = \mRQ \sR_\n ,\\
  &\sS_{\n+1/3} = \mST \sS_\n ,\\
  &\sT_{\n+1/3} = \mTU \sT_\n + \kST \sS_\n ,\\
  &\sU_{\n+1/3} = \sU_\n + (1 -  \mTU) \sT_\n + \kSU \sS_\n
  +(1- \mOU)\sO_\n +  \kPU \sP_\n 
   + \kQU \sQ_\n   + \kRU \sR_\n ,\\
  &\sV_{\n+1/3} = \sV_\n ,\\
  &\sW_{\n+1/3} = \sW_\n ,
\end{align}
where $\mult_{\j\k}=\,\e^{ -\alp{\j\k} \dt }$ and
\begin{align}
  \kPO =&  \frac{\alPO(\mPO - \mOU)  }{\alOU - \alPO} ,\\
  \kQO =& \frac{\alPO \alQP(\mQP - \mOU) }{(\alPO - \alQP)(\alOU - \alQP)}-\frac{\alPO \alQP (\mPO - \mOU)}{(\alPO - \alQP)(\alOU - \alPO)} ,\\
  \kRO =&  - \frac{\alPO \alQP \alRQ(\mQP - \mOU) }{(\alQP - \alRQ)(\alPO - \alQP)(\alOU -\alQP)} 
  + \frac{\alPO \alQP \alRQ(\mPO - \mOU) }{(\alQP - \alRQ)(\alPO - \alQP)(\alOU -\alPO)}    + \nonumber ,\\
  & +\frac{\alPO \alQP \alRQ (\mRQ - \mOU)}{(\alQP - \alRQ)(\alPO - \alRQ)(\alOU - \alRQ)}
  -\frac{\alPO \alQP \alRQ(\mPO - \mOU) }{(\alQP - \alRQ)(\alPO - \alRQ)(\alOU - \alPO)}    ,\\
  \kQP =& \frac{\alQP(\mQP - \mPO) }{\alPO - \alQP} ,\\
  \kRP =& -\frac{\alQP \alRQ(\mQP - \mPO) }{(\alQP - \alRQ)(\alPO - \alQP)} 
  + \frac{\alQP \alRQ(\mRQ - \mPO) }{(\alQP - \alRQ)(\alPO - \alRQ)} ,\\
  \kRQ =& -\frac{\alRQ(\mQP - \mRQ)}{\alQP - \alRQ}  ,\\
  \kST =& -\frac{\alST(\mTU - \mST)}{\alTU- \alST} ,\\
  \kSU =& 1 + \frac{\alST \mTU - \alTU \mST}{\alTU- \alST} ,\\
  \kPU =&  1- \frac{\alOU \mPO - \alPO \mOU}{\alOU - \alPO} ,\\
  \kQU =&  \frac{\alPO}{\alPO - \alQP} \left( 1- \frac{\alOU \mQP - \alQP \mOU}{\alOU - \alQP}\right)
  -\frac{ \alQP }{\alPO - \alQP}   \left(1 - \frac{\alOU \mPO - \alPO \mOU}{\alOU - \alPO}\right) ,\\
  \kRU =& - \frac{\alPO \alRQ }{(\alQP - \alRQ)(\alPO - \alQP)}  \left( 1- \frac{\alOU \mQP - \alQP \mOU}{\alOU -\alQP} \right)  + \nonumber \\
  & + \frac{\alQP \alRQ }{(\alQP - \alRQ)(\alPO - \alQP)}  \left( 1-  \frac{\alOU \mPO - \alPO \mOU}{\alOU -\alPO}\right)  + \nonumber \\
  & +\frac{ \alPO \alQP }{(\alQP - \alRQ)(\alPO - \alRQ)} \left( 1-  \frac{\alOU \mRQ - \alRQ \mOU}{\alOU - \alRQ} \right) -\nonumber \\
  & -\frac{ \alQP \alRQ }{(\alQP - \alRQ)(\alPO - \alRQ)}  \left( 1- \frac{\alOU \mPO - \alPO \mOU}{\alOU - \alPO} \right) .
\end{align}

At the second sub-step, the mapping \eq{sol-subs-2} is calculated by
solving the following initial-value problem, defined by the matrix $\mA_1$ \eq{A1},
\begin{align}
 &\ddt{\sO} = - \alOP \sO,		&& \sO(0)=\sO_{\n+1/3}, \eqlabel{A1-O} \\
 &\ddt{\sP} = \alOP \sO - \alPQ \sP,	&& \sP(0)=\sP_{\n+1/3}, \eqlabel{A1-P} \\
 &\ddt{\sQ} = \alPQ \sP - \alQR \sQ,	&& \sQ(0)=\sQ_{\n+1/3}, \eqlabel{A1-Q} \\
 &\ddt{\sR} = \alQR \sQ,  	    	&& \sR(0)=\sR_{\n+1/3}, \eqlabel{A1-R} \\
 &\ddt{\sS} = \alTS \sT, 		&& \sS(0)=\sS_{\n+1/3}, \eqlabel{A1-S} \\
 &\ddt{\sT} = \alUT \sU - \alTS \sT,	&& \sT(0)=\sT_{\n+1/3}, \eqlabel{A1-T} \\
 &\ddt{\sU} = - \alUT \sU,	    	&& \sU(0)=\sU_{\n+1/3}, \eqlabel{A1-U} \\
 &\ddt{\sV} = 0,	    		&& \sV(0)=\sV_{\n+1/3}, \eqlabel{A1-V} \\
 &\ddt{\sW} = 0,			&& \sW(0)=\sW_{\n+1/3}. \eqlabel{A1-W} 
\end{align}
Here we proceed similar to the first sub-step. 
We note that the equations \eq{A1-O} and \eq{A1-U} are decoupled, and
solve them to get $\sO(\t)$ and $\sU(\t)$.
The result for $\sO(\t)$ is substituted to \eq{A1-P} and the result for $\sU(\t)$ is
substituted to \eq{A1-T} to obtain closed initial value problems, which
are solved to produce $\sP(\t)$ and $\sT(\t)$.
The solution for $\sP(\t)$ is substituted to the \eq{A1-Q} to give $\sQ(\t)$.
Finally, we substitute the $\sQ(\t)$ into \eq{A1-R} and $\sT(\t)$ to \eq{A1-S}
which yield $\sR(\t)$ and $\sS(\t)$ respectively.
With the obvious solution to $\sV$ and $\sW$ the mapping is as
follows:
\begin{align}
  &\sO_{\n+2/3} = \mOP \sO_{\n+1/3}, \\
  &\sP_{\n+2/3} = \lOP \sO_{\n+1/3} + \mPQ \sP_{\n+1/3}, \\
  &\sQ_{\n+2/3} = \lOQ \sO_{\n+1/3} + \lPQ \sP_{\n+1/3} + \mQR \sQ_{\n+1/3}, \\
  &\sR_{\n+2/3} = \lOR \sO_{\n+1/3} + \lPR \sP_{\n+1/3} +  (1 -  \mQR)\sQ_{\n+1/3}  +\sR_{\n+1/3},      \\
  &\sS_{\n+2/3} = \lUS \sU_{\n+1/3} + (1- \mTS) \sT_{\n+1/3} +   \sS_{\n+1/3}, \\
  &\sT_{\n+2/3} = \lUT \sU_{\n+1/3} + \mTS \sT_{\n+1/3}, \\
  &\sU_{\n+2/3} = \mUT \sU_{\n+1/3}, \\
  &\sV_{\n+2/3} = \sV_{\n+1/3}, \\
  &\sW_{\n+2/3} = \sW_{\n+1/3},
\end{align}
where $\mult_{\j\k}=\,\e^{ -\alp{\j\k} \dt }$ again, and
\begin{align}
  &\lOP = \frac{\alOP( \mOP - \mPQ )}{\alPQ -\alOP} ,\\
  &\lOQ = \frac{\alPQ \alOP( \mOP -\mQR)}{(\alPQ -\alOP)(\alQR - \alOP)}  - \frac{\alPQ \alOP(\mPQ - \mQR)}{(\alPQ -\alOP)(\alQR - \alPQ)}  ,\\
  &\lPQ = \frac{\alPQ(\mPQ-\mQR)}{\alQR - \alPQ} ,\\
  &\lOR = 1+\frac{ \alPQ( \alOP \mQR-\alQR \mOP)}{(\alPQ -\alOP)(\alQR - \alOP)} - \frac{ \alOP(\alPQ \mQR-\alQR \mPQ)}{(\alPQ -\alOP)(\alQR - \alPQ)} ,\\
  &\lPR = 1+ \frac{\alPQ \mQR- \alQR \mPQ}{\alQR - \alPQ} ,\\
  &\lUS = 1+ \frac{\alUT \mTS-\alTS \mUT}{\alTS - \alUT} ,\\
  &\lUT = \frac{\alUT(\mUT-\mTS)}{\alTS - \alUT} .
\end{align}

The third sub-step mapping \eq{sol-subs-3} is calculated by
solving the following initial-value problem, defined by the matrix $\mA_2$ \eq{A2},
\begin{align}
  &\ddt{\sO} = \alUO \sU ,                             && \sO(0)=\sO_{\n+2/3},\eqlabel{A2-O}\\
  &\ddt{\sP} = \alUP \sU-\alPU \sP ,                     && \sP(0)=\sP_{\n+2/3},\eqlabel{A2-P}\\
  &\ddt{\sQ} = \alTQ \sT-\alQT \sQ ,                     && \sQ(0)=\sQ_{\n+2/3},\eqlabel{A2-Q}\\
  &\ddt{\sR} = \alSR \sS-\alRS \sR ,                     && \sR(0)=\sR_{\n+2/3},\eqlabel{A2-R}\\
  &\ddt{\sS} = \alRS \sR-\alSR \sS ,                     && \sS(0)=\sS_{\n+2/3},\eqlabel{A2-S}\\
  &\ddt{\sT} = \alQT \sQ-\alTQ \sT ,                     && \sT(0)=\sT_{\n+2/3},\eqlabel{A2-T}\\
  &\ddt{\sU} = \alPU \sP+\alVU \sV-(\alUP+\alUO+\alUV) \sU,&& \sU(0)=\sU_{\n+2/3},\eqlabel{A2-U}\\
  &\ddt{\sV} = \alUV \sU+\alWV \sW-(\alVU+\alVW) \sV ,     && \sV(0)=\sV_{\n+2/3},\eqlabel{A2-V}\\
  &\ddt{\sW} = \alVW \sV-\alWV \sW ,                     && \sW(0)=\sW_{\n+2/3}.\eqlabel{A2-W}
\end{align}
Unlike the previous sub-steps, this system is not solved exactly, but
its solution is approximated by the forward Euler method as follows:
\begin{align}
  &\sO_{\n+1} = \sO_{\n+2/3} + (\alUO \sU_{\n+2/3}                                   ) \dt,\\
  &\sP_{\n+1} = \sP_{\n+2/3} + (\alUP \sU_{\n+2/3}-\alPU \sP_{\n+2/3}                    ) \dt,\\
  &\sQ_{\n+1} = \sQ_{\n+2/3} + (\alTQ \sT_{\n+2/3}-\alQT \sQ_{\n+2/3}                    ) \dt,\\
  &\sR_{\n+1} = \sR_{\n+2/3} + (\alSR \sS_{\n+2/3}-\alRS \sR_{\n+2/3}                    ) \dt,\\
  &\sS_{\n+1} = \sS_{\n+2/3} + (\alRS \sR_{\n+2/3}-\alSR \sS_{\n+2/3}                    ) \dt,\\
  &\sT_{\n+1} = \sT_{\n+2/3} + (\alQT \sQ_{\n+2/3}-\alTQ \sT_{\n+2/3}                    ) \dt,\\
  &\sU_{\n+1} = \sU_{\n+2/3} + \left[\alPU \sP_{\n+2/3}+\alVU \sV_{\n+2/3}-(\alUP+\alUO+\alUV) \sU_{\n+2/3} \right] \dt,\\
  &\sV_{\n+1} = \sV_{\n+2/3} + \left[\alUV \sU_{\n+2/3}+\alWV \sW_{\n+2/3}-(\alVU+\alVW) \sV_{\n+2/3}       \right] \dt,\\
  &\sW_{\n+1} = \sW_{\n+2/3} + (\alVW \sV_{\n+2/3}-\alWV \sW_{\n+2/3}                    ) \dt.
\end{align}
This completes the definition of the hybrid method. 

\section{Details of error analysis}

As proclaimed in the main text and as we shall see below, the local
truncation errors at time step $[\t_\n,\t_\n+\dt]$ in all three methods
are given by expressions of the form
\begin{align}
  \Err(\t_\n)\dt^2+\O(\dt^3),                     \eqlabel{error-local}
\end{align}
giving the upper estimate of a global error for the interval
$\t\in[\tmin,\tmax]$ of the first order,
\begin{align}
  \sup\limits_{[\tmin,\tmax]}\norm{\vu{}^{\mathrm{exact}}-\vu{}^{\mathrm{numeric}}}
  \leq\sup\limits_{[\tmin,\tmax]}(\Err(\t)) \; (\tmax-\tmin) \dt
  + \O(\dt^2),                                    \eqlabel{error-global}
\end{align}
where the estimates of the coefficients $\Err$ are different for the three methods.

To obtain these estimates, let us consider the quasi-linear
system~\eq{markov-generic}, rewritten as
\begin{align*}
  \ddt{\vu{}} = \mA(\Vm(\t)) \vu{} = \mA(\t) \vu{} 
\end{align*}
on the interval $\t\in[\t_\n,\t_{\n+1}]$, $\t_{\n+1}=\t_\n+\dt$.
Using matrix exponential, the result
can be written in the form
\begin{align}
  \vu{}(\t_{\n+1})
  =\exp\left[ \int\limits_{\t_\n}^{\t_\n+\dt} \mA(\tf) \,\d\tf \right]\,\vu{}(\t_\n)
  \equiv \M{}(\t_\n,\dt)\,\vu{}(\t_\n)            \eqlabel{matrix-exponential}
\end{align}
The accuracy in finding $\vu{}(\t_{\n+1})$ at a given $\vu{}(\t_{\n})$
depends on accuracy of the approximation of operator $\M{}$ and on the
norm of vector $\vu{}(\t_\n)$. Since each component of $\vu{}$ is
restricted to the interval $[0,1]$ and sum of its components is fixed
to 1, we have $\norm{\vu{}}\le1$ for any choice of norm, in which any
vector, that has exactly one component equal to unity and the rest
equal to zero, is a unit vector.

Hence from this point on we focus on the approximation of the
timestep transition operator $\M{}$.

Expanding~\eq{matrix-exponential}, first the integral, then the
exponential, in the Taylor series, we have
\begin{align*}
  \M{}&
  =\exp\left[ 
    \int\limits_{\t_\n }^{\t_{\n}+\dt} \left(
      \mA(\t_\n)+\dot\mA(\t_\n)(\tf-\t_\n)+\O((\tf-\t_\n)^2)
    \right) \,\d\tf \right]
  =\exp\left[ \mA(\t_\n)\dt+\frac12\dot\mA(\t_\n)\dt^2+\O(\dt^3) \right]
  \\&
  =1 + \left[ \mA(\t_\n)\dt+\frac12\dot\mA(\t_\n)\dt^2+\O(\dt^3) \right]
     + \frac12\left[ \mA(\t_\n)\dt+\frac12\dot\mA(\t_\n)\dt^2+\O(\dt^3) \right]^2
     + \O(\dt^3)
  \\&
  =1 + \mA(\t_\n)\dt+\frac12\left[\dot\mA(\t_\n)+\mA^2(\t_\n)\right]\dt^2 + \O(\dt^3) ,
\end{align*}
where the dot designates time differentiation. 
FE approximates this operator as
\begin{align*}
  \MEF(\t_\n,\dt) = 1 + \mA(\t_\n)\dt,
\end{align*}
hence for the principal term of the norm of the error we have
\begin{align*}
  \ErrFE=\lim\limits_{\dt\to0}\norm{\MEF-\M{}}/\dt^2=
  \frac12\norm{\mA^2+\dot\mA} 
  \leq \frac12\left( \norm{\mA}^2 + \norm{\dot\mA} \right) .
\end{align*}

For the MRL, we have 
\begin{align*}
  \MMRL = \exp\left(\mA(\t_\n)\dt\right) = 1 + \mA(\t_\n)\dt + \frac12 \mA^2(\t_\n)\dt^2
  + \O(\dt^3),
\end{align*}
therefore
\begin{align*}
  \ErrMRL=\frac12\norm{\dot\mA}
  =\frac12\norm{\mA'}\abs{\dot\Vm}
\end{align*}
where the prime designates differentiation by $\Vm$. 

The errors of the three substeps of HOS are described by the above
formulas for FE (for $\mA_2$) and for MRL (for $\mA_0$, $\mA_1$), and
in addition to those, we have the error due to operator splitting. To
estimate the latter, let us compare the exact solution with the result
of the successive application of the substeps as if they were done
exactly. Let $\mB_\m=\int\limits_{\t_\n}^{\t_\n+\dt}\mA_\m(\t)\,\d\t$,
$\m=0,1,2$. Then the exact solution is
\begin{align*}
  \M{}&
  =1+(\mB_0+\mB_1+\mB_2)+\frac12(\mB_0+\mB_1+\mB_2)^2+\O(\dt^3)
  \\&
  =1+(\mB_0+\mB_1+\mB_2)+\frac12(
  \mB_0^2+\mB_1^2+\mB_2^2
  +\mB_0\mB_1+\mB_1\mB_0
  +\mB_0\mB_2+\mB_2\mB_0
  +\mB_1\mB_2+\mB_2\mB_1
  )+\O(\dt^3),
\end{align*}
and the result of the three substeps, with $\e^{\mB_0}$ applied first and
$\e^{\mB_2}$ applied last, is
\begin{align*}
  \MOS &
  =
  \left(1+\mB_2+\frac12\mB_2^2+\O(\dt^3)\right)
  \left(1+\mB_1+\frac12\mB_1^2+\O(\dt^3)\right)
  \left(1+\mB_0+\frac12\mB_0^2+\O(\dt^3)\right)
      \\&
  =
  1+\mB_2+\mB_1+\mB_0+\frac12\left(
    \mB_2^2 + \mB_1^2 + \mB_0^2
    +2\mB_2\mB_1
    +2\mB_2\mB_0
    +2\mB_1\mB_0
    \right) + \O(\dt^3) ,
\end{align*}
so

\begin{align*}
  \MOS-\M{}=\frac12\left(
  \comm{\mB_2}{\mB1}+\comm{\mB_2}{\mB_0}+\comm{\mB_1}{\mB_0}
  \right)+\O(\dt^3)
  =\frac12\left(
    \comm{\mA_2}{\mA1}+\comm{\mA_2}{\mA_0}+\comm{\mA_1}{\mA_0}
  \right)\dt^2+\O(\dt^3)
\end{align*}
where we use the standard notation for the matrix commutator,
$\comm{\matA}{\matB}\equiv\matA\matB-\matB \matA$. Finally, by the
triangle inequality (subadditivity) of a matrix norm, the upper
estimate of the error coefficient $\ErrHOS$ is given by the sum of the
error coefficients of the three constituent steps and of the operator
splitting.

\begin{figure}[tbp]
  \centering
  \includegraphics{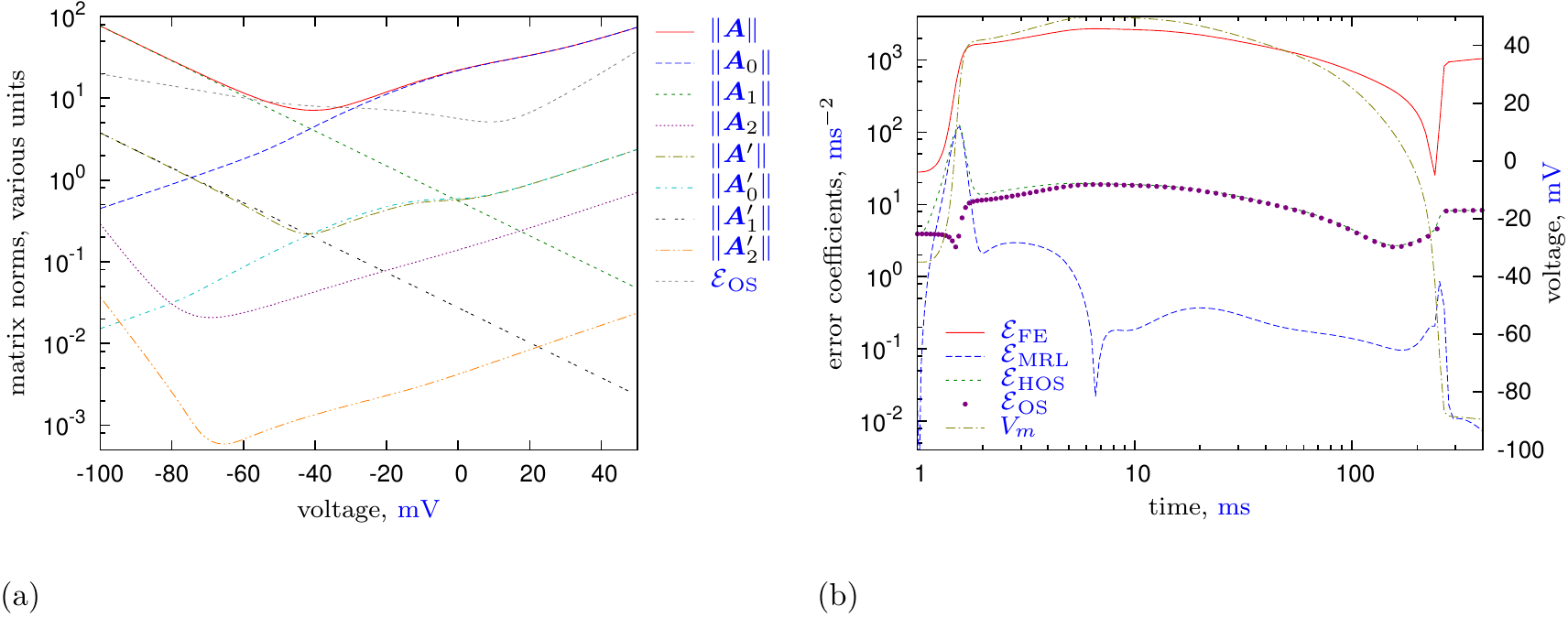}
  \caption[]{%
    (a) Matrix norms affecting the apriori error estimates, as
    functions of the transmembrane voltage. Norms
    $\norm{\mA}$,$\norm{\mA_\m}$ are in $\ms^{-1}$,
    $\norm{\mA'}$,$\norm{\mA'_\m}$,$\ErrOS$ are in $\ms^{-2}$. %
    (b) The coefficients of the apriori error estimates, as functions
    of time. $\ErrOS$ component shown by points as it overlaps with
    $\ErrHOS$ most of the time. The action potential $\Vm(\t)$ is
    shown for reference. %
  }
  \figlabel{errors}
\end{figure}

To summarize, we have the following estimates of the leading terms of
the approximation errors for the three methods as
\begin{align}
  & \ErrFE \leq \frac12\left(
    \norm{\mA}^2 + \norm{\d\mA/\d\Vm}\abs{\d\Vm/\d\t} 
  \right) \nonumber\\
  & \ErrMRL=\frac12\norm{\d\mA/\d\Vm}\abs{\d\Vm/\d\t} \nonumber\\
  & \ErrHOS \leq \frac12\abs{\d\Vm/\d\t}\left(
    \norm{\d\mA_0/\d\Vm} +
    \norm{\d\mA_1/\d\Vm} +
    \norm{\d\mA_2/\d\Vm}
  \right) + \frac12||\mA_2||^2+\ErrOS , \nonumber\\
  & \ErrOS=\frac12\norm{\comm{\mA_1}{\mA_0} +
    \comm{\mA_2}{\mA_0} + \comm{\mA_2}{\mA_1}}.   \eqlabel{err-estimates}
\end{align}

An important observation is that the apriori estimates of the errors
cannot be made based on the properties of the MC alone as they depend
on the rate of change of the voltage.

The graphs of the Frobenius norms of the matrices involved in the
estimates~\eq{err-estimates} are shown in~\fig{errors}(a). Evidently
$\norm{\mA}$ dominates other norms throughout the voltage range;
however, it is relatively small for intermediate values of $\Vm$ and
this is precisely when $\d\Vm/\d\t$ is large during a typical AP,
making the related components of the errors more significant. So a
more adequate idea of the relative magnitudes of the errors of the
three methods should take into account properties of specific
solutions. \Fig{errors}(b) shows the values of the error estimates
\eq{err-estimates} for the typical AP which was used for other
numerical illustrations in the paper. We see that the error associated
with FE is the largest of the three, with the maximal magnitude of
about $2700\,\ms^{-2}$, achieved early during the plateau of the AP,
thus guaranteeing no more than 10\% global error on a time interval of
$1\,\ms$ long for time steps as short as $\dt\approx0.04\,\us$, and its
main contributor is $\norm{\mA}^2$ rather than $\norm{\dot\mA}$. The
error associated with MRL is the smallest of the three, with the
maximal magnitude of about $118\,\ms^{-2}$, achieved during the
upstroke of the action potential, giving 10\% global accuracy on
$1\,\ms$ interval for $\dt\approx0.8\,\us$.  The error of the HOS is
intermediate between the two. Its maximum of about $125\,\ms^{-2}$,
i.e. very similar to that of MRL and achieved at the same time, as its
main contributors are the same $\dot\Vm$-dependent errors of the
exponential integrator substeps as in $\ErrMRL$. Outside the AP upstroke,
the error of HOS is dominated by the operator splitting error
$\ErrOS$, which however never exceeds $19\,\ms^{-2}$. The ratio of
the error coefficients of the two methods varies widely during the AP
solution: $\ErrFE/\ErrMRL\in(3.18,\infty)$ (remember $\ErrMRL=0$ when
$\d\Vm/\d\t=0$) and $\ErrFE/\ErrHOS\in(2.30,161)$, with the smallest
values achieved during the upstroke when the exponential solvers are
least accurate.

Clearly, the estimate of the global error given by \eq{error-global}
is over-cautious, or ``pessimistic'', as it presumes that local errors
take maximal values allowed by the matrix norms, and accumulate but
not compensate on the whole interval $[\tmin,\tmax]$. As the numerical
experiments described in the main text show, the actual errors are
much smaller. Still, the analysis done here can be useful in
identifying relative contribution of different sources of errors and
identifying ``bottlenecks''. Specifically, we see that
\begin{itemize}
\item the exponential solvers are more accurate than FE: even in the
  worst case, during the upstroke, they give two to three times
  smaller error;
\item the principal limitation of the accuracy of both exponential
  solvers is the dependence on $\dot\Vm$, which affects accuracy
  mostly during the upstroke, hence any attempts to improve the
  accuracy should in the first instance address this issue.
\end{itemize}


\end{document}